\input amstex
\documentstyle{amsppt}
\magnification=\magstep1

\loadbold
\NoBlackBoxes

\TagsAsMath
\TagsOnRight

\pageheight{8.125truein}
\pagewidth{5.5truein}
\hcorrection{0.375truein}
\vcorrection{0.5truein}

\def\ZZ{{\Bbb Z}}
\def\QQ{{\Bbb Q}}

\def\CC{{\Bbb C}}

\def\fract{\operatorname{Fract}}
\def\bfp{{\boldsymbol p}}
\def\bfq{{\boldsymbol q}}

\def\gfrak{{\frak g}}
\def\hfrak{{\frak h}}
\def\bfrak{{\frak b}}
\def\nfrak{{\frak n}}
\def\sln{\frak{sl}_n}

\def\Oq{{\Cal O}_q}
\def\Oqp{{\Cal O}_{q,p}}
\def\Mn{M_n(k)}
\def\GLn{GL_n(k)}
\def\SLn{SL_n(k)}

\def\symp{\operatorname{\frak {sp}}}
\def\eucl{\operatorname{\frak o}}

\def\Olp{\Cal O_{\lambda, \bfp}}
\def\Dlp{D_{\lambda, \bfp}}

\def\Oqkn{{\Cal O}_{\bfq}(k^n)}
\def\OlpMn{\Olp(M_n(k))}
\def\OlpMmn{\Olp(M_{m,n}(k))}
\def\OlpGLn{\Olp(GL_n(k))}
\def\OlpSLn{\Olp(SL_n(k))}
\def\OqMn{\Oq(\Mn)}

\def\OqGLn{\Oq(\GLn)}
\def\OqSLn{\Oq(\SLn)}
\def\UqgM{{U_q(\gfrak,M)}}
\def\UqbM{{U_q(\bfrak^+,M)}}

\def\UqpbM{{U_{q,p}(\bfrak^+,M)}}
\def\Uqg{{U_q(\gfrak)}}
\def\Uqb{{U_q(\bfrak^+)}}
\def\Uqn{{U_q(\nfrak^+)}}

\def\Uqpb{{U_{q,p}(\bfrak^+)}}
\def\Uqpn{{U_{q,p}(\nfrak^+)}}
\def\AnQ{A_n^{Q,\Gamma}(k)}
\def\Oqsymp{\Oq(\symp k^{2n})}
\def\Oqeucl{\Oq(\eucl k^n)}
\def\kxn{(k^\times)^n}

\def\H{{\Cal H}}

\def\E{{\Cal E}}

\def\O{{\Cal O}}

\def\SS{{\Cal S}}
\def\Ahat{A^\wedge}
\def\Hhat{\widehat{\H}}

\def\spec{\operatorname{spec}}
\def\Hspec{{\Cal H}\text{-}\roman{spec}\,}
\def\Hsp{H\text{-}\roman{spec}\,}
\def\prim{\operatorname{prim}}

\def\rann{\operatorname{r{.}ann}}
\def\lann{\operatorname{l{.}ann}}
\def\aut{\operatorname{Aut}}

\def\AlDu{{\bf 1}}
\def\AST{{\bf 2}}
\def\BrGo{{\bf 3}}
\def\qHilb{{\bf 4}}
\def\BrWa{{\bf 5}}
\def\Cald{{\bf 6}}
\def\Caldtwo{{\bf 7}}
\def\Cau{{\bf 8}}
\def\Cli{{\bf 9}}
\def\DeCK{{\bf 10}}
\def\DKPone{{\bf 11}}
\def\DeCL{{\bf 12}}
\def\DeCP{{\bf 13}}
\def\Dix{{\bf 14}}
\def\GolMi{{\bf 15}}
\def\Gdnsn{{\bf 16}}
\def\GLen{{\bf 17}}
\def\GLettwo{{\bf 18}}
\def\qaffine{{\bf 19}}
\def\specstrat{{\bf 20}}
\def\GoWa{{\bf 21}}
\def\HLevone{{\bf 22}}
\def\HLevtwo{{\bf 23}}
\def\HLT{{\bf 24}}
\def\IoMa{{\bf 25}}
\def\Jan{{\bf 26}}
\def\Jat{{\bf 27}}
\def\Jor{{\bf 28}}
\def\Jtwo{{\bf 29}}
\def\Jbook{{\bf 30}}
\def\Jfconj{{\bf 31}}
\def\KKMS{{\bf 32}}
\def\Let{{\bf 33}}
\def\LeSt{{\bf 34}}
\def\Lus{{\bf 35}}
\def\Mal{{\bf 36}}
\def\Manone{{\bf 37}}
\def\Mantwo{{\bf 38}}
\def\McRob{{\bf 39}}
\def\MoRenone{{\bf 40}}
\def\Mon{{\bf 41}}
\def\MoPa{{\bf 42}}
\def\Mus{{\bf 43}}
\def\Nor{{\bf 44}}
\def\Oh{{\bf 45}}
\def\Ohcat{{\bf 46}}
\def\OhPa{{\bf 47}}
\def\Pan{{\bf 48}}
\def\PaWa{{\bf 49}}
\def\ReVo{{\bf 50}}
\def\Ren{{\bf 51}}
\def\Res{{\bf 52}}
\def\RTF{{\bf 53}}
\def\Rin{{\bf 54}}
\def\Sas{{\bf 55}}
\def\Smi{{\bf 56}}
\def\Sud{{\bf 57}}
\def\Tak{{\bf 58}}
\def\VonII{{\bf 59}}
\def\Zha{{\bf 60}}

\topmatter

\title Prime spectra of quantized coordinate rings \endtitle

\author K. R. Goodearl \endauthor

\address Department of Mathematics, University of California, Santa
Barbara, CA 93106, USA\endaddress
\email goodearl\@math.ucsb.edu\endemail

\thanks This research was partially supported by NATO Collaborative
Research Grant 960250 and National Science
Foundation research grant DMS-9622876.\endthanks

\endtopmatter

\document

This paper is partly a report on current knowledge concerning the
structure of (generic) quantized coordinate rings and their prime
spectra, and partly propaganda in support of the conjecture that since
these algebras share many common properties, there must be a common
basis on which to treat them. The first part of the paper is expository.
We survey a number of classes of quantized coordinate rings, as well as
some related algebras that share common properties, and we record some
of the basic properties known to occur for many of these algebras,
culminating in stratifications of the prime spectra by the actions of
tori of automorphisms. As our main interest is in the generic case, we
assume various parameters are not roots of unity whenever convenient. In
the second part of the paper, which is based on
\cite{\specstrat}, we offer some support for the conjecture above, in
the form of an axiomatic basis for the observed stratifications and
their properties. At present, the existence of a suitable supply of
normal elements is taken as one of the axioms; the search for better
axioms that yield such normal elements is left as an open problem.

\head I. Quantized Coordinate Rings and Related Algebras\endhead

This part of the paper is an expository account of the prime ideal
structure of algebras on the ``quantized coordinate ring'' side of the
theory of quantum groups -- quantizations of the coordinate rings of
affine spaces, matrices, semisimple groups, symplectic and Euclidean
spaces, as well as a few related algebras -- quantized enveloping
algebras of Borel and nilpotent subalgebras of semisimple Lie algebras,
and quantized Weyl algebras. These algebras occur widely
throughout the quantum groups literature, and different papers often
investigate different versions. Thus, we begin by giving definitions of
the most general versions of which we are aware; in quoting results from
the literature, we will specify which version is under consideration.
The reader should bear in mind that many authors, when studying one
version of a quantized coordinate algebra, use results proved for
slightly different versions, on the understanding that the proofs carry
over. This is especially prevalent with regard to quantized coordinate
rings of semisimple groups; a detailed development covering the most
general case would be a welcome addition to the literature.

Fix a base field $k$ throughout. It need not be algebraically closed,
and may have arbitrary characteristic.

\head 1. Descriptions\endhead

This section is designed to be a reference source for descriptions of the
main classes of the standard quantized coordinate rings currently
studied, as well as some related algebras with similar properties. (The
reader lacking a strong stomach for generators and relations
should just skim this section and refer back to it as necessary.)
For the sake of uniformity, and to emphasize that these algebras are
deformations of classical coordinate rings, we label all quantized
coordinate rings using notations of the form
$\O_\bullet(\Cal C)$, where
$\bullet$ records one or more parameters and $\Cal C$ records the name
of the classical object. Thus, $\Oq(k^n)$ refers to a one-parameter
quantization of the coordinate ring of affine $n$-space over a field
$k$, while $\Olp(GL_n(k))$ refers to a multiparameter quantization of
the coordinate ring of $GL_n(k)$, and so on. Many different labels are
used throughout the literature for these algebras, and we do not attempt
to list the alternates here.

\definition{1.1. Quantum affine spaces} These are meant to be viewed as
deformations of the coordinate rings $\O(k^n)$. Recall that $\O(k^n)$ is
a (commutative) polynomial ring $k[x_1,\dots,x_n]$, where $x_i$ is the
$i$-th coordinate function on $k^n$. In the present case, one deforms
$\O(k^n)$ in a rather straightforward manner, by altering the
commutativity relations $x_ix_j=x_jx_i$ to ``commutativity up to
scalars'', that is, $x_ix_j= q_{ij}x_jx_i$. In order to prevent
degeneracy, one needs some assumptions on the $q_{ij}$. In particular,
they should be nonzero, and $q_{ji}$ should equal $q_{ij}^{-1}$ in order
to prevent the pair of relations $x_ix_j= q_{ij}x_jx_i$ and $x_jx_i=
q_{ji}x_ix_j$ from implying $x_ix_j=0$. Thus, one assumes that $\bfq=
(q_{ij})$ is a {\it multiplicatively antisymmetric\/}
$n{\times} n$ matrix over $k$, that is,
$q_{ii}=1$ and $q_{ji}= q_{ij}^{-1}$ for all
$i,j$. Given $\bfq$, the corresponding {\it multiparameter coordinate
ring of quantum affine
$n$-space over
$k$\/} is the $k$-al\-ge\-bra $\O_\bfq(k^n)$ generated by elements
$x_1,\dots,x_n$ subject only to the relations
$x_ix_j= q_{ij}x_jx_i$ for $i,j= 1,\dots,n$. The standard single
parameter version occurs when the $q_{ij}$ for $i<j$ are all equal to a
fixed nonzero scalar $q$; in this case, we denote the algebra
$\Oq(k^n)$. 

Quantum affine spaces occur already in the work of Manin \cite{\Manone,
Section 1, \S2 and Section 4, \S5}; a superalgebra version is discussed
in \cite{\Mantwo, \S1.4}. 
\enddefinition

\definition{1.2. Quantum matrices } The set $\Mn$ of $n\times n$
matrices over $k$, viewed as an algebraic variety, is just affine
$n^2$-space, and its coordinate ring $\O(\Mn)$ is a polynomial ring in
$n^2$ indeterminates. Hence, one might expect to deform this algebra
exactly as in (1.1). However, there is more structure in this case, and
one seeks to preserve it as far as possible. Namely, $\O(\Mn)$ is a
bialgebra with a comultiplication $\Delta : \O(\Mn)\rightarrow
\O(\Mn)\otimes \O(\Mn)$ which is effectively the transpose of the map
$\Mn\times\Mn
\rightarrow
\Mn$ given by matrix
multiplication: If $X_{ij}$ denotes the $i,j$-th coordinate function on $\Mn$, then
$\Delta(X_{ij})= \sum_{l=1}^n X_{il}\otimes X_{lj}$. In particular, one
would like deformations of $\O(\Mn)$, with appropriate sets of
generators $X_{ij}$, which are also bialgebras in which the
comultiplication of the $X_{ij}$ is given by the equation above.

In addition, multiplication of matrices with row or column vectors from
$k^n$ induces on $\O(k^n)$ structures of left and right comodule over
$\O(\Mn)$, and one would like some quantum affine spaces $\O_\bfq(k^n)$
to be left and right comodules over the deformation of $\O(\Mn)$, with
comodule structure maps behaving as in the classical case. In the single
parameter case, there is a deformation $\Oq(\Mn)$ over which $\Oq(k^n)$
becomes both a left and a right comodule in the desired manner (see
\cite{\Mantwo} and \cite{\Res}). However, if one tries to make a
multiparameter quantum affine space $\O_\bfq(k^n)$ 
into a left \underbar{and} right comodule over the desired type of
deformation of
$\O(\Mn)$, the latter deformation will usually be degenerate. To obtain
nondegenerate deformations, one must use different multiparameter
quantum affine spaces  as left and right comodules, as discovered in
\cite{\AST} and \cite{\Sud}, to which we refer the reader for a more
complete discussion. The resulting algebras can be described as follows.

Let $\bfp =
(p_{ij})$ be a multiplicatively antisymmetric $n\times n$ matrix over
$k$, and let
$\lambda$ be a nonzero element of $k$ not equal to $-1$. The
corresponding {\it multiparameter coordinate ring of quantum $n\times n$
matrices over
$k$\/} is the
$k$-al\-ge\-bra $\OlpMn$  generated by elements
$X_{ij}$ (for $i,j=1,\dots,n$) subject only to the following relations:
$$X_{\ell m}X_{ij} = \cases p_{\ell i}p_{jm}X_{ij}X_{\ell m} +
(\lambda -1)p_{\ell i}X_{im}X_{\ell j}&\quad (\ell >i,\
m>j)\\ 
\lambda p_{\ell i}p_{jm}X_{ij}X_{\ell m}&\quad (\ell >i,\ m\le j)\\ 
p_{jm}X_{ij}X_{\ell m}&\quad (\ell
=i,\ m>j).\endcases$$
(We assume $\lambda\ne0$ to avoid obvious degeneracies in the second
relation above. The assumption $\lambda\ne-1$ is needed to ensure that
$\OlpMn$ has the same Hilbert series as a commutative polynomial algebra
in $n^2$ indeterminates \cite{\AST, Theorem 1}.) The standard single
parameter version, denoted
$\OqMn$, occurs when the $p_{ij}$ for $i>j$ are all equal to a fixed
nonzero scalar $q$, and
$\lambda= q^{-2}$. In some references, such as \cite{\PaWa} and
\cite{\Smi}, the roles of $q$ and $q^{-1}$ are interchanged, and/or $q$
is squared.

Quantized coordinate rings for rectangular matrices are defined as
subalgebras of those for square matrices. For instance, if $m<n$ then
$\OlpMmn$ is defined to be the $k$-subalgebra of $\OlpMn$ generated by
those $X_{ij}$ with $i\le m$. There is a $k$-al\-ge\-bra retraction of
$\OlpMn$ onto $\OlpMmn$ whose kernel is generated by the $X_{ij}$ with
$i>m$; thus $\OlpMmn\cong \OlpMn/ \langle X_{ij} \mid i>m \rangle$.
Hence, results for $\OlpMn$ easily extend to the rectangular case; we
shall so extend results from the literature without explicit mention.

\enddefinition

\definition{1.3. Quantum general linear groups} Within the algebraic
variety $\Mn$, the set $\GLn$ of invertible matrices forms an open
subvariety, the complement of the variety defined by the vanishing of
the determinant function $D$. The coordinate ring $\O(\GLn)$
thus has the form $\O(\Mn)[D^{-1}]$. To construct a quantum analog, one
inverts a ``quantum determinant'', call it $\Dlp$, in $\OlpMn$. In order
for the inversion process to work smoothly and avoid degeneracies, the
powers of $\Dlp$ should form an Ore set. This occurs because the 
chosen $\Dlp$ turns out to be a normal element.
(Recall that a {\it normal element\/} in a ring $R$ is an
element $c$ such that $cR=Rc$.) 

As in the classical case, the quantum determinant $\Dlp$ can be defined
as a linear combination of products $X_{1,\pi(1)} X_{2,\pi(2)} \cdots
X_{n,\pi(n)}$ as $\pi$ runs through the symmetric group $S_n$. However,
the coefficients are taken to be certain products of elements $-p_{ij}$
rather than powers of
$-1$. Specifically,
$$\Dlp= \sum_{\pi\in S_n} \biggl( \prod\Sb 1\le i<j\le n\\ \pi(i)>\pi(j)
\endSb  (-p_{\pi(i),\pi(j)}) \biggr) X_{1,\pi(1)} X_{2,\pi(2)} \cdots
X_{n,\pi(n)}.$$ 
The motivation for this choice of $\Dlp$ comes from the construction of
a quantized version of the exterior algebra on $k^n$; the $n$-th
``quantum exterior power'' of $k^n$ is 1-di\-men\-sional, spanned by $\Dlp$
(see \cite{\AST} for details). In the single parameter case, the quantum
determinant $D_q$ in $\OqMn$ can be expressed as
$$D_q= \sum_{\pi\in S_n} (-q)^{\ell(\pi)} X_{1,\pi(1)} X_{2,\pi(2)}
\cdots X_{n,\pi(n)},$$
where $\ell(\pi)$ denotes the {\it length\/} of the permutation $\pi$,
that is, the minimum length of an expression for $\pi$ as a product of
adjacent transpositions $(i,i+1)$. (Cf\. \cite{\PaWa, Lemma 4.1.1} but
interchange $q$ and $q^{-1}$.)

It has been computed that $\Dlp$ is a normal element of $\OlpMn$
\cite{\AST, Theorem 3}; in fact,
$$\Dlp X_{ij}= \lambda^{j-i} \biggl( \prod_{l=1}^n p_{jl}p_{li} \biggr)
X_{ij}\Dlp$$
for all $i,j$. In the single parameter case, $D_q$ is  central
(e.g., \cite{\PaWa, Theorem 4.6.1}). It follows that the set of
nonnegative powers of
$\Dlp$ is a right and left Ore set in $\OlpMn$. The {\it multiparameter
coordinate ring of quantum $\GLn$\/} is now defined to be the
localization $\Olp(GL_n(k))= \OlpMn[\Dlp^{-1}]$; in the single parameter
case, $\OqMn[D_q^{-1}]$ is denoted $\OqGLn$.
\enddefinition

\definition{1.4. Quantum special linear groups} Note that the set
$\SLn$ of $n\times n$ matrices with determinant 1 forms a closed
subvariety of $\Mn$, defined by the single equation $D=1$. The
coordinate ring $\O(\SLn)$ thus has the form $\O(\Mn)/(D-1)$. To
construct a quantum analog, we would thus factor out from $\OlpMn$ the
ideal generated by $\Dlp-1$. However, unless $\Dlp$ is central, such a
factor would be degenerate. Hence, quantum special linear groups are
only defined for special choices of the parameters.

From the normality relations for $\Dlp$, we see that $\Dlp$ is central
if and only if
$$\lambda^i \prod_{l=1}^n p_{il}= \lambda^j \prod_{l=1}^n p_{jl}$$
for all $i,j$. In this case, the {\it multiparameter
coordinate ring of quantum $\SLn$\/} is defined to be the factor algebra
$$\Olp(SL_n(k))= \OlpMn/ \langle \Dlp-1 \rangle.$$
In the single parameter
case, where
$D_q$ is automatically central, the factor $\OqMn/ \langle D_q-1
\rangle$ is denoted $\OqSLn$.
\enddefinition

\definition{1.5. Quantum semisimple groups} A systematic development of
quantized coordinate rings for simple algebraic groups corresponding to
the  Dyn\-kin diagrams $A_n$, $B_n$, $C_n$, $D_n$ was given by
Reshetikhin, Takhtadzhyan, and Fadeev in \cite{\RTF}. These algebras were
expressed in terms of the entries of certain ``$R$-matrices'' attached
to quantizations of the corresponding simple Lie algebras. Calculations
of generators and relations for these quantized coordinate rings can be
found in, e.g.,
\cite{\Tak}. Quantized coordinate rings for the exceptional groups,
however, were not computed, due to the lack of explicit $R$-matrices in
those cases (cf\. \cite{\RTF, Remark 14, p\. 212}; the $G_2$ case is
considered in \cite{\Sas}). The now standard approach to these algebras
involves a quantization of the classical duality between the coordinate
ring of a semisimple algebraic group and the enveloping algebra of its
Lie algebra. Thus, one constructs quantized enveloping algebras of
semisimple Lie algebras and {\it defines\/} the quantized coordinate
rings of the corresponding semisimple groups as certain Hopf algebra
duals. We outline  this procedure.

{\bf (a)} Let $\gfrak$ be a complex
semisimple Lie algebra, and let $G$ be a connected complex semisimple
algebraic group with Lie algebra $\gfrak$. Both $\gfrak$ and $G$ play
symbolic roles; they serve mainly as suggestive labels. The only
data needed from these objects are: a (symmetrized) Cartan matrix
$(d_ia_{ij})$ for $\gfrak$; root and weight lattices $Q\subseteq P$
corresponding to choices of Cartan subalgebra $\hfrak$ and root system
for $\gfrak$; a lattice $L$ lying between $Q$ and $P$
corresponding to the character group of a maximal torus of $G$ with Lie
algebra $\hfrak$; choices of simple roots $\alpha_1,\dots,\alpha_n$ and
fundamental weights $\omega_1,\dots,\omega_n$; and the (unique) bilinear
pairing $(-,-):
P\times Q \rightarrow \ZZ$ such that $(\omega_i,\alpha_j)=
\delta_{ij}d_i$ for all $i,j$. There is a unique extension of $(-,-)$
to a symmetric bilinear form $P\times P\rightarrow \QQ$, which we also
denote by
$(-,-)$. 

Different authors base quantized coordinate algebras of $G$ on different
quantizations of the enveloping algebra $U(\gfrak)$. Thus we first
describe a general one-parameter quantization of $U(\gfrak)$, as in
\cite{\DKPone, \S0.3}. Let $M$ be a lattice lying between $Q$ and $P$
(independent of $L$). Then let $q\in k$ be a nonzero scalar, set
$q_i=q^{d_i}$ for $i=1,\dots,n$, and assume that these $q_i\ne\pm1$. (In
much of the literature on quantized enveloping algebras,
$k$ is taken to be $\CC$, or a rational function field $\CC(q)$ or
$\QQ(q)$, or the algebraic closure of one of these rational function
fields.)  Let $U^0$ be a copy of
the group algebra $kM$, written as the $k$-al\-ge\-bra with basis $\{
k_\lambda \mid \lambda\in M\}$ where $k_0=1$ and $k_\lambda k_\mu=
k_{\lambda+\mu}$ for $\lambda,\mu\in M$. The {\it single parameter
quantized enveloping algebra of $\gfrak$} associated with the above
choices is the $k$-al\-ge\-bra $\UqgM$ generated by $U^0$ and elements
$e_1,\dots,e_n$, $f_1,\dots,f_n$ satisfying the following relations for
$\lambda\in M$ and  $i,j=1,\dots,n$:
$$\allowdisplaybreaks\gather k_\lambda e_ik_\lambda^{-1} =
q^{(\lambda,\alpha_i)} e_i
\qquad\text{and}\qquad k_\lambda f_ik_\lambda^{-1} =
q^{-(\lambda,\alpha_i)} f_i\\
e_if_j-f_je_i= \delta_{ij} \frac{k_{\alpha_i}-k_{\alpha_i}^{-1}}
{q_i-q_i^{-1}}\\
\sum_{l=1}^{1-a_{ij}} (-1)^l {{1-a_{ij}}\choose l}_{q_i}
e_i^{1-a_{ij}-l} e_je_i^l =0 \qquad(i\ne j)\\
\sum_{l=1}^{1-a_{ij}} (-1)^l {{1-a_{ij}}\choose l}_{q_i}
f_i^{1-a_{ij}-l} f_jf_i^l =0 \qquad(i\ne j),\endgather$$
where ${{1-a_{ij}}\choose l}_{q_i}$ is a $q_i$-binomial coefficient. The
most commonly studied case is $U_q(\gfrak,Q)$, denoted just $\Uqg$. This
algebra is generated by $e_i,f_i,k_i$ for $i=1,\dots,n$, where $k_i=
k_{\alpha_i}$. At the other extreme, $U_q(\gfrak,P)$ is
often denoted $\check U_q(\gfrak)$.

The algebra $\UqgM$ is in fact a Hopf algebra, with comultiplication
$\Delta$, counit $\epsilon$, and antipode $S$ such that
$$\xalignat3 \Delta(k_\lambda) &= k_\lambda \otimes k_\lambda
&\epsilon(k_\lambda) &=1 &S(k_\lambda) &= k_\lambda^{-1}\\
\Delta(e_i) &= e_i\otimes 1 +k_{\alpha_i}\otimes e_i
&\epsilon(e_i) &=0 &S(e_i) &= -k_{\alpha_i}^{-1}e_i\\
\Delta(f_i) &= f_i\otimes k_{\alpha_i}^{-1} +1\otimes f_i
&\epsilon(f_i) &=0 &S(f_i) &= -f_ik_{\alpha_i}.\endxalignat$$
Hence, one can define
the {\it finite dual\/} $\UqgM^\circ$ as in, e.g., \cite{\Mon, Definition
1.2.3}; this is a $k$-linear subspace of $\UqgM^*$ which becomes a Hopf
algebra using the transposes of the multiplication, comultiplication, and
antipode of
$\UqgM$ (e.g., \cite{\Mon, Theorem 9.1.3}).

{\bf (b)} Single-parameter quantizations of $\O(G)$ are defined as
subalgebras of $\UqgM^\circ$ generated by the ``coordinate functions''
of certain ``highest weight'' $\UqgM$-modules. It turns out that the
resulting algebras are independent of the choice of $M$. Hence, we first
describe quantizations of $\O(G)$ as subalgebras of $\Uqg^\circ$; this
has the advantage that no roots of $q$ are required. For descriptions in
terms of $\UqgM^\circ$, see part (c). To avoid problems with certain
calculations, one assumes that $\text{char}(k)\ne2$, and also
$\text{char}(k)\ne3$ in case $\gfrak$ has a component of type $G_2$.

Here we only give a quantization of $\O(G)$ for the case that $q$ is not
a root of unity. The root of unity case requires defining a suitable
algebra over a Laurent polynomial ring $k_0[t^{\pm1}]$ and then
specializing $t$ to $q$ (see, e.g., \cite{\Lus, Sections 7,8} or
\cite{\DeCL, Section 4}).

For each
$\lambda\in P^+$, there is a finite dimensional simple
$\Uqg$-module $V(\lambda)$ with {\it highest weight\/} $\lambda$, that
is,
$V(\lambda)$ is generated by an element $u_\lambda$ satisfying $k_\mu
u_\lambda= q^{(\mu,\lambda)}u_\lambda$ for all $\mu\in Q$ and
$e_iu_\lambda =0$ for $i=1,\dots,n$ (see, e.g.,
\cite{\Jan, Theorem 5.10}). For
$v\in V(\lambda)$ and $f\in V(\lambda)^*$, let $c^{V(\lambda)}_{f,v}\in
\Uqg^\circ$ denote the {\it coordinate function\/} defined by the rule
$c^{V(\lambda)}_{f,v}(u)= f(uv)$. The {\it single parameter quantized
coordinate ring of $G$} is the $k$-subalgebra $\Oq(G,L)$ of $\Uqg^\circ$
generated by the $c^{V(\lambda)}_{f,v}$ for $\lambda\in L^+$, $f\in
V(\lambda)^*$, and
$v\in V(\lambda)$. There is some redundancy in the notation $\Oq(G,L)$,
since $L$ is determined by $G$, but we prefer to emphasize $L$ in this
way. (Thus, the pair $(G,L)$ is used as a label for that connected
semisimple group with Lie algebra $\gfrak$ and weight lattice $L$.) The
most commonly studied case corresponds to a {\it simply connected\/}
group. This is the case $L=P$, and in this case we write $\Oq(G)$ for
$\Oq(G,P)$.

When $G=SL_n(\CC)$, one obtains the algebra $\Oq(M_n(\CC))/\langle
D_q-1\rangle$ as in (1.4) (e.g.,
\cite{\HLevone, Theorem 1.4.1}; replace
$q^2$ by $q$). Thus, there is no ambiguity in the notation
$\Oq(SL_n(\CC))$.

{\bf (c)} In order to exhibit $\Oq(G,L)$ as a subalgebra of
$\UqgM^\circ$, one needs sufficient roots of $q$ in $k$ so that
$q^{(M,L)}\subset k$. For $\lambda\in L^+$, there is then a simple finite
dimensional $\UqgM$-module $V(M,\lambda)$ with highest weight $\lambda$,
where now the highest weight vector $u_\lambda$ must satisfy $k_\mu
u_\lambda= q^{(\mu,\lambda)}u_\lambda$ for all $\mu\in M$. Let us
 write $\Oq(G,L,M)$ for the $k$-subalgebra of $\UqgM^\circ$
generated by the coordinate functions of the $V(M,\lambda)$ for
$\lambda\in L^+$. Each $V(M,\lambda)$ is also simple as a $\Uqg$-module;
thus $V(M,\lambda)$ becomes $V(\lambda)$ by restriction, and coordinate
functions of $V(M,\lambda)$ map to coordinate functions of $V(\lambda)$
by restriction. Therefore the restriction map $\UqgM^\circ \rightarrow
\Uqg^\circ$ induces a $k$-al\-ge\-bra homomorphism $\Oq(G,L,M) \rightarrow
\Oq(G,L)$. This map is an isomorphism, because
$$\Oq(G,L,M)= \bigoplus_{\lambda\in L^+} C^{V(M,\lambda)}
\qquad\text{and}\qquad \Oq(G,L)= \bigoplus_{\lambda\in L^+}
C^{V(\lambda)}$$
where $C^{V(M,\lambda)}$ and $C^{V(\lambda)}$ denote the $k$-linear
spans of the coordinate functions on $V(M,\lambda)$ and $V(\lambda)$,
respectively (cf\. \cite{\Jtwo, \S2.2}, \cite{\Jbook, \S1.4.13},
\cite{\HLT,
\S3.3}). One also needs the fact that $C^{V(M,\lambda)}$ and
$C^{V(\lambda)}$ have the same $k$-di\-men\-sion, since both are isomorphic
to $V(\lambda)\otimes V(\lambda)^*$.

For the reasons sketched above, we can -- and do -- choose to work with
$\Oq(G,L)$. In \cite{\HLT, \S3.3}, for instance, the corresponding
algebra -- there denoted $\CC_q[G]$ -- is defined as $\Oq(G,L,L)$.

{\bf (d)} Multiparameter versions of $\Oq(G,L)$ are obtained by twisting
the multiplication. Let $p:L\times L \rightarrow k^\times$ be an {\it
alternating bicharacter\/}, that is, 
$$\align p(\lambda,\lambda) &=1\\
p(\mu,\lambda) &= p(\lambda,\mu)^{-1}\\
p(\lambda,\mu+\mu') &= p(\lambda,\mu) p(\lambda,\mu') \endalign$$
for $\lambda,\mu,\mu' \in  L$. There is a natural
$( L\times L)$-bigrading on $\Oq(G,L)$ (cf\. \cite{\HLT,
\S3.3}), and  the multiplication on $\Oq(G,L)$ can be twisted via
$p$ as in \cite{\HLT, \S2.1}; the new multiplication, call it
$\cdot$, is determined by the rule
$$a\cdot b= p(\lambda,\lambda')
p(\mu,\mu')^{-1} ab \qquad \text{for\ } a\in \Oq(G,L)_{\lambda,\mu}
\text{\ and\ } b\in \Oq(G,L)_{\lambda',\mu'}.$$
The vector space $\Oq(G,L)$, equipped with this new multiplication, is
called a {\it multiparameter quantized coordinate ring of $G$}, denoted
$\Oqp(G,L)$. In the case $L=P$, we  write simply $\Oqp(G)$.
\enddefinition

\definition{1.6. Quantized enveloping algebras of Borel and nilpotent
subalgebras} Each of the quantized enveloping algebras $\UqgM$ contains
subalgebras which are viewed as quantizations of the enveloping algebras
of Borel or nilpotent subalgebras of $\gfrak$. Unlike $\UqgM$, these
algebras behave much like $\Oq(G)$, and they play prominent roles in the
study of $\Oq(G)$. Hence, we include them here under the rubric of
``algebras similar to quantized coordinate rings'', along with quantized
Weyl algebras (see the following subsection). 

Carry over the objects and  notation from (1.5a), in particular the Lie
algebra
$\gfrak$, lattices $Q\subseteq M\subseteq P$, the powers
$q_i=q^{d_i}\ne\pm1$, and the group algebra $U^0\subseteq \UqgM$. Let
$\bfrak^+$ and
$\nfrak^+$ denote the positive Borel and nilpotent subalgebras of
$\gfrak$ corresponding to the choices above.  The {\it single parameter
quantized enveloping algebra of $\bfrak^+$} associated with the above
choices is the $k$-al\-ge\-bra $\UqbM$ generated by $U^0$ and 
$e_1,\dots,e_n$. The
corresponding {\it single parameter
quantized enveloping algebra of $\nfrak^+$} is the $k$-subalgebra $\Uqn$
of $\UqbM$ generated by $e_1,\dots,e_n$; this algebra is often denoted
$\Uqg^+$ (it is, of course, independent of the choice of
$M$).

There is an $(M\times M)$-bigrading on the algebra $\UqbM$ such that
$k_\lambda
\in \UqbM_{-\lambda,\lambda}$ for $\lambda\in M$ and $e_i\in
\UqbM_{-\alpha_i,0}$ for $i=1,\dots,n$ (see \cite{\HLT, Corollary 3.3}).
Obviously, $\Uqn$ is a bigraded subalgebra. Given an alternating
bicharacter
$p: M\times M
\rightarrow k^\times$, we can twist the multiplications in $\UqbM$ and
$\Uqn$ via $p$ just as for $\Oq(G,L)$ above. The resulting algebras are
the {\it multiparameter quantized enveloping algebras of $\bfrak^+$ and
$\nfrak^+$}, denoted $\UqpbM$ and $\Uqpn$. As in (1.5a), we omit $M$ from
the notation in the case $M=Q$.

There are two  cases in which a quantized Borel is known to be a
homomorphic image of a quantum semisimple group. Namely, when $q$ is
transcendental over $\QQ$, the algebra
$U_q(\bfrak^+,P)$ is isomorphic to a factor algebra of
$\Oq(G)$ \cite{\Jbook, Corollary 9.2.12}, and for $q\in\CC$ not a root
of unity, $\UqpbM$ is a homomorphic image of
$\O_{q,p^{-1}}(G,M)^{\text{op}}$ \cite{\HLT, Proposition 4.6}.
\enddefinition

\definition{1.7. Quantized Weyl algebras} Recall that the Weyl algebra
$A_n(k)$ can (at least in characteristic zero) be viewed as the algebra
of polynomial differential operators on affine $n$-space, i.e., the
algebra of all linear partial differential operators with polynomial
coefficients on the coordinate ring $\O(k^n)$. Quantized versions of
$A_n(k)$ arose in Maltsiniotis' development of ``quantum differential
calculus''
\cite{\Mal}, and can be viewed as algebras of ``partial $q$-difference
operators'' on quantum affine spaces.

Let $Q= (q_1,\dots, q_n)$ be a vector in $(k^\times)^n$, and let $\Gamma
=(\gamma_{ij})$ be a multiplicatively antisymmetric $n\times n$ matrix
over $k$. The {\it multiparameter quantized Weyl algebra of degree $n$
over
$k$\/} is the $k$-al\-ge\-bra $\AnQ$ generated by
elements $x_1,y_1,\dots, x_n,y_n$ subject only to the following
relations:
$$\alignat2 y_iy_j &=\gamma_{ij}y_jy_i &\qquad &(\text{all\ } i,j)\\
x_ix_j &=q_i\gamma_{ij}x_jx_i &\qquad &(i<j)\\ 
x_iy_j &=\gamma_{ji}y_jx_i &\qquad &(i<j)\\ 
x_iy_j &=q_j\gamma_{ji}y_jx_i &\qquad &(i>j)\\ 
x_jy_j &=1+q_jy_jx_j
+\sum_{l<j} (q_l-1)y_lx_l &\qquad &(\text{all\ }j). \endalignat$$
See \cite{\Jor, \S2.9} for a discussion of how to view elements of $\AnQ$
as $q$-dif\-fer\-ence operators on $\O_\Gamma(k^n)$.
\enddefinition

\definition{1.8. Quantum symplectic spaces} The classical linear,
symplectic, and orthogonal groups all act on the same affine space,
$k^n$. We have already mentioned above that the relations in quantum
matrices are partly chosen so that quantum affine spaces become
comodules over quantum matrices, in a ``standard'' manner. It follows
that quantum affine spaces are also comodules over quantum special
linear groups. In order to develop an analogous situation for quantum
symplectic and  orthogonal groups, one needs, as it turns out,
different deformations of $\O(k^n)$ than quantum affine spaces. These
new algebras are called quantum symplectic and Euclidean spaces; they
have (to our knowledge) only been  defined in single parameter versions
to date.

Let $q$ be a nonzero scalar in $k$ and $n$ a positive integer. Set
$$\align (1,2,\dots,n,n',(n-1)',\dots,1') &= (1,2,\dots,n,n+1,\dots,
2n)\\
(\rho(1),\rho(2), \dots,\rho(2n)) &= (n,n-1,\dots,1, -1,-2, \dots,-n)\\
(\epsilon_1,\dots,\epsilon_n,\epsilon_{n+1}, \dots,\epsilon_{2n}) &=
(1,\dots,1,-1,\dots,-1). \endalign$$
The one-parameter {\it coordinate
ring of quantum symplectic
$2n$-space over
$k$\/} (cf\. \cite{\RTF, Definition
14} or \cite{\Mus, \S1.1}) is the $k$-al\-ge\-bra $\Oqsymp$
generated by elements
$x_1,\dots,x_{2n}$ satisfying the following relations:
$$\alignat2 x_ix_j &= qx_jx_i &\qquad &(i<j;\ j\ne i')\\
x_ix_{i'} &= x_{i'}x_i+ (1-q^2) \sum_{l=1}^{i'-1} q^{\rho(i')-\rho(l)}
\epsilon_{i'}
\epsilon_l x_lx_{l'}  &\qquad &(i\le n) \endalignat$$
(cf\. \cite{\RTF, p\. 210}). A simpler set of relations was
found by Musson \cite{\Mus, \S1.1} (cf\.
\cite{\Oh, \S1.1}):
$$\alignat2 x_ix_j &= qx_jx_i &\qquad &(i<j;\ j\ne i')\\
x_ix_{i'} &= q^2x_{i'}x_i +(q^2-1)\sum_{l=1}^{i-1}
q^{l-i}x_lx_{l'} &\qquad &(i\le n). \endalignat$$
\enddefinition

\definition{1.9. Quantum Euclidean spaces} Let $n\ge2$ be an integer
and
$q$  a nonzero scalar in
$k$, such that $q$ has a square root in $k$ in case $n$ is odd. Set
$m=\lfloor n/2
\rfloor$, the integer part of $n/2$, and set $i'=n+1-i$ for $1\le i\le
n$. Further, set
$$\multline (\rho(1),\dots,\rho(n))\\
= {\cases (m-\tfrac12, m-\tfrac32, \dots,
\tfrac12, 0, -\tfrac12, -\tfrac32, \dots, -m+\tfrac12) &\qquad (n\text{\
odd})\\
(m-1,m-2, \dots,1,0,0,-1,-2,\dots, -m+1) &\qquad (n\text{\ even}).
\endcases} \endmultline$$
The one-parameter {\it coordinate
ring of quantum Euclidean
$n$-space over
$k$\/} (cf\. \cite{\RTF, Definition
12} or \cite{\Mus, \S2.1}) is the $k$-al\-ge\-bra $\Oqeucl$
generated by elements
$x_1,\dots,x_n$ satisfying the following relations:
$$\align 
x_ix_j &= qx_jx_i \tag"$(i<j;\ j\ne i')$\quad"\\
x_ix_{i'} &= x_{i'}x_i+ \dfrac{1-q^2}{1+q^{n-2}} \biggl( q^{n-2}
\sum_{l=1}^{i'-1} q^{\rho(i')-\rho(l)} x_lx_{l'} -\sum_{l=i'}^n
q^{\rho(i')-\rho(l)} x_lx_{l'} \biggr)\\
 & \tag"$(i\le m)$.\quad" \endalign$$
A simpler set of relations was given in \cite{\Mus, \S\S2.1, 2.2}:
$$\align 
x_ix_j &= qx_jx_i \tag"$(i<j;\ j\ne i')$\quad"\\
x_ix_{i'} &= x_{i'}x_i+ (1-q^2) \sum_{l=i+1}^m q^{l-i-2} x_lx_{l'}
\tag"$(i\le m;\ n\text{\ even})$\quad"\\
x_ix_{i'} &= x_{i'}x_i+ (1-q)q^{m-i-(1/2)}x_{m+1}^2+ (1-q^2)
\sum_{l=i+1}^m q^{l-i-2} x_lx_{l'}\\
 & \tag"$(i\le m;\ n\text{\ odd})$.\quad" \endalign$$
For an alternative presentation of $\Oqeucl$, see \cite{\Ohcat, Example 5;
\OhPa, Definition 3.1}.
\enddefinition

\head 2. Some common properties\endhead

We record some fundamental properties common to the algebras discussed
in Section 1. The cases we list are those for which we have located
references in which a result appears, or from which it follows
readily. Many cases which we conjecture should be included must be
omitted at present, because they have not (to our knowledge) been
addressed in the literature. We leave it to the reader to identify
appropriate ``missing'' cases, and to formulate the corresponding
conjectures. Perhaps more important than filling in specific cases is
the  problem of developing general theorems which could verify
properties for all these algebras simultaneously.

\definition{2.1. Noetherian domains}
\enddefinition

\proclaim{Theorem} All of the following algebras are noetherian domains:
$\Oqkn$; $\OlpMmn$; $\OlpGLn$; $\OlpSLn$; $\Oqp(G,L)$ for $q$ transcendental
over $\QQ$ or $q\in\CC$ not a root of unity; $\Uqpb$; $\Uqpn$; $\AnQ$;
$\Oqsymp$; $\Oqeucl$.
\endproclaim

\demo{Proof} In the cases $\Oqkn$; $\OlpMmn$; $\AnQ$; $\Oqsymp$;
$\Oqeucl$, one has only to observe that the algebra is an iterated skew
polynomial ring over $k$. This is clear for $\Oqkn$; for the other
cases, see \cite{\AST, pp\. 890-891}, \cite{\Jor, \S\S2.1, 2.8},
\cite{\Mus, \S\S1.2, 2.3}. It follows that the localization $\OlpGLn$ is
a noetherian domain, and that the factor algebra $\OlpSLn$ is
noetherian. That the latter is a domain is proved in \cite{\LeSt,
Corollary}.

For $\Oq(G)$ in the stated cases, see \cite{\Jtwo, Lemma 3.1, Proposition
4.1} and
\cite{\Jbook, Lemma 9.1.9, Proposition 9.2.2}. (An alternate proof of noetherianity
is given in \cite{\qHilb, Corollary 5.6}.) The desired properties of
$\Oq(G,L)$ are proved by the same arguments, and they carry over to
$\Oqp(G,L)$ by graded ring methods \cite{\HLT, Remark, p\. 80} or by
twisting results \cite{\Zha, Propositions 5.1, 5.2}.

By \cite{\DeCK, Corollary 1.8}, $\Uqg$ is a domain, and so $\Uqb$ and
$\Uqn$ are domains. There are several ways to see that $\Uqb$ and $\Uqn$
are noetherian. For instance, it follows from the existence of a PBW
basis for $\Uqg$ (e.g.,
\cite{\Jan, Theorem 4.21}) that $\Uqg$ is free as a right $\Uqb$-module
and as a right
$\Uqn$-module. Hence, the poset of left ideals of either subalgebra
embeds in the poset of left ideals of $\Uqg$, via $I\mapsto \Uqg I$.
This shows that $\Uqb$ and $\Uqn$ are left noetherian. There exists an
antiautomorphism $\tau$ on $\Uqg$ such that $\tau(e_i)=e_i$ for 
$i=1,\dots,n$ and $\tau(k_\lambda)= k_\lambda^{-1}$ for $\lambda\in Q$
(e.g., \cite{\Jan, Lemma 4.6}). Since $\tau$ obviously stabilizes $\Uqb$
and
$\Uqn$, these subalgebras must be right as well as left noetherian.

On the other hand, in view of \cite{\DeCK, Proposition 1.7}, given
either $A_0=\Uqb$ or $A_0=\Uqn$, there is a sequence of algebras
$A_0,A_1,\dots,A_N$ such that each $A_{i+1}$ is the associated graded
ring of $A_i$ with respect to some $\ZZ^+$-filtration, and such that
$A_N$ is an iterated skew polynomial ring over either $U^0$ or $k$.
Since $A_N$ is noetherian, it follows from standard filtered/graded
techniques (e.g., \cite{\McRob, Theorem 1.6.9}) that $A_0$ is noetherian.

Finally, the twisting results mentioned above (\cite{\Zha, Propositions
5.1, 5.2}) imply that $\Uqpb$ and $\Uqpn$ are noetherian domains.
\qed\enddemo

\definition{2.2. Complete primeness of prime factors} For any set
$\{\alpha_i\}$ of nonzero scalars in $k$, let $\langle \alpha_i \rangle$
denote the multiplicative subgroup of $k^\times$ generated by the
$\alpha_i$.
\enddefinition

\proclaim{Theorem} Let $A$ be one of the following algebras: 
$\Oqkn$ with $\langle q_{ij}\rangle$ torsionfree; $\OlpMmn$ or $\OlpGLn$ or
$\OlpSLn$ with $\langle
\lambda,\, p_{ij}\rangle$ torsionfree; $\Oq(G)$ or $\Uqb$ with $q\in\CC$ not a
root of unity; $\Uqn$ with $q$ transcendental over $\QQ$; $\AnQ$ with $\langle
q_i,\,
\gamma_{ij}\rangle$ torsionfree; $\Oqsymp$ or $\Oqeucl$ with $q$ not a
root of unity.
 Then all prime
ideals of
$A$ are completely prime, i.e., all prime factor rings of $A$ are
domains.
\endproclaim

\demo{Proof} The cases $\Oqkn$; $\OlpMmn$; $\AnQ$; $\Oqsymp$; $\Oqeucl$
follow from a general result about prime ideals of certain iterated skew
polynomial rings: \cite{\GLettwo, Theorem 2.3}. See \cite{\GLettwo,
Theorem 2.1}, \cite{\Gdnsn, Theorem 5.1}, \cite{\Cau, Proposition
II.1.2}, \cite{\Mus, Corollary 1.2, Theorem 2.3} for further details.
The cases $\OlpGLn$ and $\OlpSLn$ follow immediately. In addition,
$\Uqn$ is an iterated $q$-skew polynomial ring over $k$ (see \cite{\Rin,
Section 5} for the case $k=\QQ(q)$ and extend scalars), and so in this
case also the result follows from \cite{\GLettwo, Theorem 2.3} -- see
\cite{\Rin, Corollary to Theorem 3}. Finally, the cases $\Oq(G)$ and
$\Uqb$ are proved in \cite{\Jtwo, Theorems 11.4, 11.5}. \qed\enddemo

\definition{2.3. Division rings of fractions}
\enddefinition

\proclaim{Theorem} Let $A$ be one of the algebras $\OlpMmn$;
$U_q(\sln(k))^+$ with $q$ not a root of unity; $\Oq(G)$ or $U_q(\bfrak^+,P)$ or
$\Uqn$ with $q\in\CC$ not a root of unity;
$\AnQ$;
$\Oqsymp$;
$\Oqeucl$. Then $\fract A\cong \fract \O_\bfq(K^t)$ for some $\bfq$, some
field $K\supseteq k$, and some $t$.

Now let $A$ be either $\OqMn$ with $q$  not a root of unity, or $\AnQ$
with $\langle q_i,\,
\gamma_{ij}\rangle$ torsionfree. If $P$ is any prime ideal of $A$, then
$\fract(A/P) \cong \fract \O_\bfq(K^t)$ for some $\bfq$, some
field $K\supseteq k$, and some $t$. \endproclaim

\demo{Proof} For the first set of cases, see \cite{\MoPa, Theorem 1.24},
\cite{\AlDu, Th\'eor\`eme 2.15}, \cite{\IoMa, Theorem 3.5}, \cite{\Caldtwo,
Theorems 3.3, 3.1},
\cite{\Jor, \S3.1},
\cite{\Mus, Theorems 1.3, 2.4}. The case of $\fract \OqMn= \fract \OqGLn$
was treated earlier in \cite{\Cli, Proposition 5} and \cite{\Pan, Theorem
3.8}. The case of $\Uqn$ for $q$ transcendental over $\CC$ was also obtained in
\cite{\Jfconj}. For the cases of prime factors of
$\OqMn$ and
$\AnQ$, see
\cite{\Cau, Th\'eor\`emes II.2.1, III.3.2.1}. \qed\enddemo

\definition{2.4. Catenarity} Recall that (the prime spectrum of) a ring
$A$ is {\it catenary\/} if for any comparable prime ideals $P\supset Q$
in $A$, all saturated chains of prime ideals from $P$ to $Q$ have the
same length. Our discussion of this property will involve the following
concept. The ring $A$, or its prime spectrum
$\spec A$, is said to have {\it normal separation} provided the
following condition holds: For any proper inclusion $P\supset Q$ of
prime ideals of $A$, the factor $P/Q$ contains a nonzero element which
is normal in $A/Q$. In other words, there must exist an element $c\in
P\setminus Q$ such that $c+Q$ is normal in $A/Q$; we refer to the latter
condition by saying that $c$ is {\it normal modulo $Q$}.
\enddefinition

\proclaim{Theorem} The following algebras are catenary: $\Oqkn$;
$\Oq(GL_n(\CC))$ and $\Oq(SL_n(\CC))$ with $q\in\CC$ not a root of
unity; $\Uqn$ with $q$ transcendental over $\QQ$; $\AnQ$ with no $q_i$ a
root of unity; $\Oqsymp$ and $\Oqeucl$ with $q$ not a root of unity.
\endproclaim

\demo{Proof} These cases all follow from a general theorem based on
Gabber's methods: If $A$ is an affine, noetherian, Auslander-Gorenstein,
Cohen-Mac\-aul\-ay algebra with finite Gelfand-Kirillov dimension, and if
$\spec A$ is normally separated, then $A$ is catenary \cite{\GLen,
Theorem 1.6}. See \cite{\GLen, Theorems 2.6, 3.13, 4.5, 4.8} and
\cite{\Ohcat, Corollaries 12, 13} for the individual cases listed.
\qed\enddemo

\definition{2.5. Normal separation} As indicated in (2.4), normal
separation plays a key role in proving catenarity. Normal separation in
a ring $A$ also implies Jategaonkar's strong second layer condition
(cf\. \cite{\Jat, Theorem 3.3.16, Proposition 8.1.7, discussion p\.
225}),  from which it follows, in particular, that every module has a
finite filtration such that associated primes in adjacent layers are
linked in $\spec A$ (cf\. \cite{\Jat, Theorem 9.1.2}). Other influences
of the second layer condition on the representation theory of $A$ are
discussed in \cite{\BrWa} and \cite{\Jat}. Finally, we mention that
normal separation, in conjunction with several other properties typical
of quantum coordinate rings, leads to a description of links and
cliques in terms of lattices of automorphisms \cite{\BrGo, Section 3}.
\enddefinition

\proclaim{Theorem} The following algebras have normal separation:
$\Oqkn$; $\Oq(G)$ with $q$
transcendental over $\QQ$\/; $\Oqp(G,L)$ with $q\in\CC$ not a root of
unity; $U_q(\bfrak^+,P)$ and $\Uqn$ over a rational function field
$\CC(q)$;
$U_q(\bfrak^+,P)$ with
$q$ transcendental over $\QQ$\/; $\UqpbM$ with
$q\in\CC$ not a root of unity;
$\AnQ$ with no $q_i$ a
root of unity; $\Oqsymp$ and $\Oqeucl$ with $q$ not a root of unity.
\endproclaim

\demo{Proof} For $\Oqkn$, $\AnQ$, $\Oqsymp$ and $\Oqeucl$, see
\cite{\GLen, Corollary 2.4, Theorem 3.12} and \cite{\Ohcat, Theorem 10}.
The algebras $\Uqn$ and $U_q(\bfrak^+,P)$ over $\CC(q)$ are actually {\it
polynormal\/}, i.e., every ideal has a normalizing sequence of
generators \cite{\Cald, Corollaires 3.2, 3.3}.

Normal separation for $A=\Oq(G)$ or $A=\Oqp(G,L)$ is a consequence of
results in \cite{\Jbook} and \cite{\HLT}, as follows (cf\. \cite{\BrGo,
Theorem 5.8}). In these cases, there is a partition 
$$\spec A= \bigsqcup_{w\in W\times W} \spec_w A,$$
where $W$ is the Weyl group of $G$, described in \cite{\Jbook, Corollary
9.3.9} and \cite{\HLT, Corollary 4.5}. For each $w$, there are an ideal
$I_w$ of $A$ and an Ore set $\E_w$ of nonzero normal elements in $A/I_w$
such that
$$\spec_w A= \{P\in \spec A\mid P\supseteq I_w\ \text{and\ } (P/I_w)\cap
\E_w =\varnothing\}$$
\cite{\Jbook, Proposition 10.3.2}, \cite{\HLT, Theorem 4.4}.
In particular, normal separation follows for primes $P\supset Q$ such
that $Q\in\spec_w A$ but $P\notin \spec_w A$. Further, the
localization $(A/I_w)[\E_w^{-1}]$ has {\it central separation\/}: for
any primes $P'\supset Q'$, there exists a central element in
$P'\setminus Q'$. In one of our cases, this follows from a sequence of
results in \cite{\Jbook}, as shown in \cite{\BrGo, \S5.7}; in the other
case, it is known that all ideals of $(A/I_w)[\E_w^{-1}]$ are centrally
generated \cite{\HLT, Theorem 4.15}. Central separation in
$(A/I_w)[\E_w^{-1}]$ together with normality of the elements of $\E_w$
yields normal separation for primes $P\supset Q$ in the same $\spec_w A$
\cite{\BrGo, Proposition 1.7}. An alternate proof of normal separation
for $\Oq(G)$, based on Hopf algebra technology, is given in \cite{\Let,
Proposition 2.4}.

Finally, the second $U_q(\bfrak^+,P)$ case and the $\UqpbM$ case follow
from the cases of $\Oq(G)$ and $\Oqp(G,L)$ (see the last paragraph of
(1.6)).
\qed\enddemo

\definition{2.6. The Dixmier-Moeglin equivalence} Recall that a prime
ideal $P$ of a noetherian $k$-al\-ge\-bra $A$ is {\it rational\/} provided
the center of $\fract(A/P)$ is algebraic over $k$. Recall also that $P$
is {\it locally closed\/} in $\spec A$ provided the singleton $\{P\}$ is
closed in some Zariski-neighborhood, i.e., the intersection of the primes
properly containing $P$ is larger than $P$. One says that the algebra
$A$ satisfies the {\it Dixmier-Moeglin equivalence\/} if the sets of
rational prime ideals, locally closed prime ideals, and primitive ideals
all coincide (cf\. \cite{\Ren}). The advantage of this equivalence is
that when it holds, the primitive ideals of $A$ can be identified without
constructing any irreducible representations.
\enddefinition

\proclaim{Theorem} The Dixmier-Moeglin equivalence holds in the
following algebras: $\Oqkn$; $\OlpMmn$; $\OlpGLn$; $\OlpSLn$; $\Oq(G)$
with $q$ transcendental over $\QQ$ and $k$ algebraically closed;
$\Oqp(G,L)$ with $q\in\CC$ not a root of unity; $U_q(\bfrak^+,P)$ with
$q$ transcendental over $\QQ$; $\UqpbM$ with
$q\in\CC$ not a root of unity; $\AnQ$ with the $q_i$
not roots of unity; $\Oqsymp$ and $\Oqeucl$ with $q$ not a root of
unity. \endproclaim

\demo{Proof} For $\Oqkn$, $\OlpMmn$, $\AnQ$, $\Oqsymp$ and $\Oqeucl$, see
\cite{\qaffine, Corollary 2.5, Theorem 3.2} and \cite{\specstrat,
Theorems 5.3, 5.5, 5.8, 5.11}. The cases $\OlpGLn$ and $\OlpSLn$ follow
directly. The case  $\Oqp(G,L)$ follows from the results of
\cite{\HLT} exactly as in the case $\Oq(SL_n(\CC))$, which was given
explicitly in \cite{\HLevtwo, Theorem 4.2}. The case $\Oq(G)$
follows from results in \cite{\Jbook}, as noted in \cite{\specstrat,
\S2.4}. Finally, the cases $U_q(\bfrak^+,P)$ and $\UqpbM$ follow
from the previous cases.
\qed\enddemo

\head 3. Stratified spectra\endhead

In this section, we give a more extensive discussion of another feature
common to many quantized coordinate algebras -- a stratification of the
prime spectrum in which each stratum is a classical scheme, in fact the
prime spectrum of a (commutative) Laurent polynomial algebra. Such a
stratification was first discovered in $\spec \Oq(SL_3(\CC))$ by Hodges
and Levasseur \cite{\HLevone}; they soon extended it to $\Oq(SL_n(\CC))$
\cite{\HLevtwo}. This was generalized by Joseph to $\Oq(G)$ \cite{\Jtwo,
\Jbook}, and finally extended to $\Oqp(G,L)$ by Hodges, Levasseur, and
Toro \cite{\HLT}. In all these cases, the result followed from a long
sequence of involved calculations specific to the algebra at hand.
Later, it was noticed that some features of these stratifications could
be tied to the action of a torus of automorphisms (cf\. \cite{\BrGo,
Section 5 and Proposition 1.9}), and similar stratifications were
observed in other quantized algebras (e.g., \cite{\BrGo, Section 4},
\cite{\qaffine, \S2.1, Theorem 2.3}). A general development of this type
of stratification was begun in \cite{\specstrat}.

We begin with an outline of the main features of the stratification
of $\spec \Oqp(G,L)$, followed by a discussion of the axiomatic framework
into which this stratification nicely fits. Then we exhibit appropriate
tori of automorphisms for the algebras from Section 1, and indicate the
cases in which the details of the stratification have been worked out.

\definition{3.1} Let $A$ be either the algebra $\Oq(G)$ with $q$
transcendental over $\QQ$, or $\Oqp(G,L)$ with $q\in\CC$ not a root of
unity. As noted in (2.5), there is a partition 
$$\spec A= \bigsqcup_{w\in W\times W} \spec_w A,$$
where $W$ is the Weyl group of $G$, described in \cite{\Jbook, Corollary
9.3.9} and \cite{\HLT, Corollary 4.5}. For each $w$, there are an ideal
$I_w$ of $A$ and an Ore set $\E_w$ of nonzero normal elements in $A/I_w$
such that
$$\spec_w A= \{P\in \spec A\mid P\supseteq I_w\ \text{and\ } (P/I_w)\cap
\E_w =\varnothing\}.$$
In fact, $\E_w$ is finitely generated, and so if $e_w$ is the product
(in some order) of a finite set of generators for $\E_w$, we can write
$$\spec_wA= \{P\in\spec A\mid P\supseteq I_w\} \cap \{P\in\spec A\mid
e_w\notin P\},$$
an intersection of a closed and an open set. Thus $\spec_w A$ is a
locally closed subset of $\spec A$. Let $A_w=(A/I_w)[\E_w^{-1}]$. In the
case of
$\Oqp(G,L)$ over $\CC$, the
center $Z(A_w)$ is a Laurent polynomial ring over $k$
 \cite{\HLT, Theorem 4.14}, the ideals of $A_w$ are centrally generated
\cite{\HLT, Theorem 4.15}, and it follows that $\spec_w A$ is
homeomorphic to $\spec Z(A_w)$. That such a homeomorphism also exists in
the case of $\Oq(G)$ over $k\supseteq \QQ(q)$ is not proved in
\cite{\Jbook} (a Laurent polynomial subalgebra $Z_w\subseteq Z(A_w)$ is
studied instead \cite{\Jbook, \S10.3.3}), but this can be deduced from
the results there together with what we shall prove in Part II.

If $\prim A$ denotes the {\it primitive spectrum\/} of $A$ (the set of
primitive ideals, equipped with the Zariski topology), we obtain a
corresponding partition
$$\prim A= \bigsqcup_{w\in W\times W} \prim_w A,$$
where each $\prim_w A= \prim A\cap \spec_w A$ is a locally closed subset
of $\prim A$. Further, $\prim_w A$ is precisely the set of maximal
elements of $\spec_w A$ (this is given explicitly in \cite{\Jbook,
Theorem 10.3.7} for the first case, and is implicit in \cite{\HLT,
Section 4} for the second; see also \cite{\BrGo, Corollary 1.5}). Over
an algebraically closed base field, a maximal torus $H$ of $G$ acts
naturally as automorphisms of $A$ \cite{\Jbook, \S10.3.8}, \cite{\HLT,
\S3.3}, and the sets $\prim_w A$ are precisely the $H$-orbits in $\prim
A$ \cite{\Jbook, Theorem 10.3.8}, \cite{\HLT, Theorem 4.16}.
\enddefinition

\definition{3.2} The picture of $\spec A$ outlined in (3.1) can be
conveniently organized in terms of the action of the group $H$. There is
a unique minimal element $J_w$ in $\spec_w A$ for each $w$
(\cite{\Jbook, Proposition 10.3.5}; implicit in \cite{\HLT, Section 4}),
and it is easily checked that $J_w$ is $H$-stable (cf\. \cite{\BrGo,
\S5.4}). Over an algebraically closed base field, it then follows that
the $J_w$ are precisely the $H$-stable prime ideals of $A$, and that
each $\spec_w A$ consists precisely of those prime ideals $P$ such that
the intersection of the $H$-orbit of $P$ equals $J_w$ \cite{\BrGo,
Proposition 1.9}.

These observations provide a framework that can be set up with respect
to any group of automorphisms of any ring, as follows. As we shall see
later, the special properties of this framework in the cases discussed
in (3.1) appear in many other algebras as well.
\enddefinition

\definition{3.3. $H$-strat\-i\-fi\-ca\-tions} Let $H$ be a group acting as
automorphisms on a ring $A$. Recall that an {\it $H$-prime ideal\/} of
$A$ is any proper $H$-stable ideal $J$ of $A$ such that a product of
$H$-stable ideals is contained in $J$ only when one of the factors is
contained in $J$. We shall write $\Hsp A$ for the set of $H$-prime
ideals of $A$. For any ideal $I$ of $A$, let $(I:H)$ denote the
intersection of the $H$-orbit of $I$, that is,
$$(I:H)= \bigcap_{h\in H} h(I).$$
Alternatively, $(I:H)$ is the largest $H$-stable ideal of $A$ contained
in
$I$.

Observe that if $P$ is a prime ideal of $A$, then $(P:H)$ is an
$H$-prime ideal. On the other hand, if $A$ is noetherian and $I$ is an
$H$-prime ideal, it is easily checked that $I$ is semiprime and that the
prime ideals minimal over $I$ form a single $H$-orbit (cf\.
\cite{\GolMi, Remarks 4*,5*, p\. 338}). In this case, $I=(P:H)$ for any
prime $P$ minimal over $I$.

The {\it $H$-stratum\/} of any $H$-prime ideal $J$ is the set
$$\spec_J A= \{P\in \spec A\mid (P:H)=J\}.$$
One defines $H$-strata in any subset of $\spec A$ by intersection with
these strata. In particular, the $H$-stratum of $J$ in $\prim A$ is the
set
$$\prim_J A= \prim A\cap \spec_J A.$$
Since $(P:H)\in \Hsp A$ for each $P\in\spec A$, the $H$-strata partition
$\spec A$, that is,
$$\spec A= \bigsqcup_{J\in\Hsp A} \spec_J A.$$
We refer to this partition as {\it the $H$-strat\-i\-fi\-ca\-tion of $\spec A$}.
Similarly, the partition of $\prim A$ into $H$-strata $\prim_J A$ is
called the $H$-strat\-i\-fi\-ca\-tion of $\prim A$.
\enddefinition

\definition{3.4} Finite $H$-strat\-i\-fi\-ca\-tions enjoy the minimal
topological properties required of stratifications in algebraic and
differential geometry, as follows. A {\it stratification\/} of an
algebraic variety $X$ is defined in
\cite{\KKMS, p\. 56} to be a partition of $X$ into a disjoint union of
finitely many locally closed subsets (called {\it strata\/}) such that the
closure of each stratum is a union of strata.
\enddefinition

\proclaim{Lemma} Let $H$ be a group acting as automorphisms on a ring
$A$, and assume that $\Hsp A$ is finite.

{\rm (a)} Each $H$-stratum in $\spec A$ is locally closed.

{\rm (b)} The closure of each $H$-stratum in $\spec A$ is a union of
$H$-strata.

{\rm (c)} For each nonnegative integer $d$, the union of the $H$-strata
corresponding to the elements of height at most $d$ in $\Hsp A$ is open
in $\spec A$. 
\endproclaim

\demo{Proof} (a) Let $J\in\Hsp A$, and let $J'$ be the intersection of
the $H$-primes properly containing $J$ (take $J'=A$ if $J$ is maximal in
$\Hsp A$). Since $J$ is $H$-prime and $\Hsp A$ is finite, $J'\supsetneq
J$. Thus
$$\spec_J A= \{P\in\spec A\mid P\supseteq J\} \cap \{P\in \spec A\mid
P\not\supseteq J'\}.$$

(b) If $J\in\Hsp A$, the closure of $\spec_J A$ is just the set $V(J)$ of
those prime ideals of $A$ containing $J$, and $V(J)$ equals the union of the
$H$-strata $\spec_K A$ as $K$ runs over those $H$-primes containing $J$.
 
(c) Let $Y$ be the set of elements of $\Hsp A$ of height greater than
$d$, that is, those $H$-primes $K$ from which there 
descends a chain $K=K_0 \supsetneq K_1\supsetneq \cdots\supsetneq
K_{d+1}$ in $\Hsp A$. The union of the $H$-strata $\spec_JA$ for $J$ of
height at most $d$ in $\Hsp A$ is the set
$$\{P\in\spec A\mid P\not\supseteq \cap Y\},$$
which is open in $\spec A$. \qed\enddemo

\definition{3.5. Tori of automorphisms} Whether or not our base field
$k$ is algebraically closed, we shall refer to any finite direct product
of copies of the multiplicative group $k^\times$ as an {\it (algebraic)
torus\/}. For each of the algebras $A$ discussed in Section 1, there is
a naturally occurring torus $\H$ acting as $k$-al\-ge\-bra automorphisms on
$A$, as follows. The action of an element $h\in\H$ on an element $a\in
A$ will be denoted $h{.}a$; we leave it to the reader to verify in each
case that there exist $k$-al\-ge\-bra automorphisms $h{.}(-)$ as described,
and that the rule $h\mapsto h{.}(-)$ is a group homomorphism from $\H$
to $\aut A$.

Many of these tori arise as groups of `winding automorphisms' of a
bialgebra or Hopf algebra, as follows. Suppose that $A$ is a bialgebra
over $k$, and let $\Ahat$ be the set of linear characters on $A$, that is, $k$-al\-ge\-bra
homomorphisms $A\rightarrow k$. This set forms a monoid under
convolution, the counit of $A$ being the identity element of $\Ahat$. For
$\chi\in\Ahat$, there are $k$-al\-ge\-bra endomorphisms
$\theta^l_\chi$ and $\theta^r_\chi$ on $A$ given by the rules
$$\theta^l_\chi(a)= \sum_{(a)} a_1\chi(a_2) \qquad\text{and}\qquad
\theta^r_\chi(a)= \sum_{(a)} \chi(a_1)a_2.$$
When $\chi$ is invertible in $\Ahat$, say with inverse $\chi'$, the
endomorphisms
$\theta^l_\chi$ and $\theta^r_\chi$ are automorphisms of $A$, with
inverses
$\theta^l_{\chi'}$ and $\theta^r_{\chi'}$ respectively, and we refer to
them as the left and right {\it winding automorphisms\/} of $A$
associated with $\chi$. The maps
$\chi\mapsto \theta^l_\chi$ and $\chi\mapsto \theta^r_\chi$ are a
homomorphism and an anti-homomorphism, respectively, from the group of
units of $\Ahat$ to $\aut A$. The facts outlined above are well known
when $A$ is a Hopf algebra, in which case $\Ahat$ is a group under
convolution (see, e.g., \cite{\Jbook, \S1.3.4}).

If $C$ is a right comodule algebra over the bialgebra $A$ (e.g., see
\cite{\Mon, Definition 4.1.2}), with structure map $C\rightarrow
C\otimes A$ written as $c\mapsto \sum_{(c)} c_0\otimes c_1$, then each
(invertible) character $\chi\in\Ahat$ induces a $k$-al\-ge\-bra endomorphism
(automorphism) of $C$ according to the rule $c\mapsto \sum_{(c)}
c_0\chi(c_1)$. The analogous statements hold for a left comodule algebra
over $A$.

{\bf 1. Quantum affine spaces.} For $A= \Oqkn$, take $\H=\kxn$ and 
$$(\alpha_1,\dots,\alpha_n){.}x_i= \alpha_ix_i$$
for $i=1,\dots,n$. These automorphisms can be viewed as winding
automorphisms arising when $A$ is made into a right comodule algebra
over the bialgebra $\O_{\lambda,\bfq}(M_n(k))$, or into a left comodule
algebra over $\O_{\lambda,\lambda^{-1}\bfq}(M_n(k))$.

{\bf 2. Quantum matrices.} For the algebra $A= \OlpMmn$, take $\H=
(k^\times)^m\times\kxn$ and
$$(\alpha_1,\dots,\alpha_m, \beta_1,\dots,\beta_n){.}X_{ij}= \alpha_i
\beta_j X_{ij}$$
for $i=1,\dots,m$ and $j=1,\dots,n$. In the case $m=n$, the
automorphisms above are compositions of left and right winding
automorphisms. Namely, for $\gamma= (\gamma_1,\dots,\gamma_n)\in \kxn$
there is a character $\chi_\gamma\in \Ahat$ such that
$\chi_\gamma(X_{ij})= \delta_{ij}\gamma_iX_{ij}$ for all $i,j$, and the
automorphism displayed above can be expressed as $\theta^l_{\chi_\beta}
\theta^r_{\chi_\alpha}$ where $\beta=(\beta_1,\dots,\beta_n)$ and
$\alpha= (\alpha_1,\dots,\alpha_n)$. The map $\gamma\mapsto \chi_\gamma$
embeds the torus $\kxn$ in the group of units of $\Ahat$; it is an
isomorphism when $\lambda$ and $\bfp$ are sufficiently generic.

{\bf 3. Quantum general linear groups.} For $A=\OlpGLn$, we can take the
 torus $\H= \kxn\times\kxn$ acting on $\OlpMn$ as above, and extend the
given automorphisms from
$\OlpMn$ to $A$, because $\Dlp$ is an $\H$-ei\-gen\-vec\-tor.

{\bf 4. Quantum special linear groups.} For $A=\OlpSLn$ (assuming a
suitable choice of $\lambda, \bfp$ so that $\Dlp$ is central), we must
take a subgroup of the torus used for $\OlpMn$, namely the stabilizer
of $\Dlp$:
$$\H= \{ (\alpha_1,\dots,\alpha_n, \beta_1,\dots,\beta_n)\in
\kxn\times\kxn \mid \alpha_1\alpha_2\cdots\alpha_n
\beta_1\beta_2\cdots\beta_n= 1\}.$$
The automorphisms of $\OlpMn$ corresponding to elements of $\H$ all fix
$\Dlp-1$ and therefore induce automorphisms of $A$.

{\bf 5. Quantum semisimple groups.} Let $A=\Oqp(G,L)$ as described in
(1.5). In the case $k=\CC$, one works with a maximal torus of $G$ with
Lie algebra $\hfrak$ and character group $L$ \cite{\HLT, \S3.1}. In the
general case, it is simplest to take the set $\H$ of group homomorphisms
$L\rightarrow k^\times$; this is a group under pointwise multiplication,
and it is a torus because $L$ is a lattice. An action of $\H$ on $A$
derives from an embedding of $\H$ into the character group $\Ahat$.
Namely, there is a monic group homomorphism $\iota : \H\rightarrow \Ahat$
such that
$$\iota(h)(c^{V(\lambda)}_{f,v})= h(\mu)f(v)$$
for $\lambda\in P^+$, $\mu\in P$, $v\in V(\lambda)_\mu$, $f\in
V(\lambda)^*$ (see
\cite{\HLT, \S3.3}).

Therefore there is a group homomorphism $\H\rightarrow \aut A$ such that
$h\mapsto \theta^l_{\iota(h)}$, which provides an action of $\H$ on
$A$. This is the action used in \cite{\Jbook} and \cite{\HLT}. For some
purposes, however, an action of $\H\times\H$ is needed, given by
$(h_1,h_2)\mapsto \theta^l_{\iota(h_1)} \theta^r_{\iota(h_2)}$ (see
\cite{\BrGo, \S5.2}).

{\bf 6. Quantized enveloping algebras of Borel and nilpotent
subalgebras.} For $A=\Uqpn$, take $\H=\kxn$ and 
$$(\alpha_1,\dots,\alpha_n){.}e_i= \alpha_ie_i$$
for $i=1,\dots,n$. For $A=\UqpbM$, let $\H_M$ be the group of group
homomorphisms $M\rightarrow k^\times$ under pointwise multiplication
and take $\H=\H_M\times\kxn$, with
$$\align (h,\alpha_1,\dots,\alpha_n){.}k_\mu &= h(\mu)k_\mu\\
(h,\alpha_1,\dots,\alpha_n){.}e_i &= \alpha_ie_i \endalign$$
for $\mu\in M$ and $i=1,\dots,n$.

{\bf 7. Quantized Weyl algebras.} For $A=\AnQ$, take $\H=\kxn$ and
$$\align (\alpha_1,\dots,\alpha_n){.}x_i &= \alpha_ix_i\\
(\alpha_1,\dots,\alpha_n){.}y_i &= \alpha_i^{-1}y_i \endalign$$
for $i=1,\dots,n$.

{\bf 8. Quantum symplectic spaces.} For $A=\Oqsymp$, take $\H=\kxn$ and
$$\align (\alpha_1,\dots,\alpha_n){.}x_i &= \alpha_ix_i\\
(\alpha_1,\dots,\alpha_n){.}x_{i'} &= \alpha_i^{-1}x_{i'} \endalign$$
for $i=1,\dots,n$.

{\bf 9. Quantum Euclidean spaces.} For $A=\Oqeucl$, set $m= \lfloor
n/2\rfloor$ and take $\H=(k^\times)^m$, with
$$\alignat2 (\alpha_1,\dots,\alpha_m){.}x_i &= \alpha_ix_i &\qquad
&(i=1,\dots,m)\\ 
(\alpha_1,\dots,\alpha_m){.}x_{i'} &= \alpha_i^{-1}x_{i'} &\qquad
&(i=1,\dots,m)\\ 
(\alpha_1,\dots,\alpha_m){.}x_{m+1} &= x_{m+1} &\qquad &(\text{if\ } n\
\text{is odd}) \endalignat$$
\enddefinition

\definition{3.6} We close this part of the paper by recording cases in
the literature where the prime spectrum of an algebra $A$ from Section
1, stratified by a torus $\H$ from (3.5), fits into the pattern
discussed in (3.1) and (3.2). Although the results quoted from
\cite{\specstrat} do not provide Ore sets of \underbar{normal} elements,
they can be modified to do so in all but one of the cases of interest
here, as we shall show in Part II.
\enddefinition

\proclaim{Theorem} Let $A$ be one of the following algebras: $\Oqkn$;
$\OlpMmn$ or $\OlpGLn$ or $\OlpSLn$ with $\lambda$ not a root of unity;
$\Oq(G)$ with
$q$ transcendental over
$\QQ$ and $k$ algebraically closed;
$\Oqp(G,L)$ with $q\in\CC$ not a root of unity; $\AnQ$ with the $q_i$
not roots of unity; $\Oqsymp$ or $\Oqeucl$ with $q$ not a root of
unity. Let $\H$ be the corresponding torus acting on $A$ as in {\rm
(3.5)}. Then $\Hspec A$ is finite, and for each $J\in \Hspec A$ the
following hold:

{\rm (a)} There exists an Ore set $\E_J$ of regular elements in $A/J$
such that $\spec_J A$ is homeomorphic to $\spec Z((A/J)[\E_J^{-1}])$ via
localization and contraction.

{\rm (b)} $Z((A/J)[\E_J^{-1}])$ is a commutative Laurent
polynomial algebra over an extension field of $k$.

{\rm (c)} $\prim_J A$ equals the set of maximal elements of $\spec_J A$.

{\rm (d)} If $k$ is algebraically closed, then $\prim_J A$ consists of a
single $\H$-orbit.
\endproclaim

\demo{Proof} For $\Oqkn$, see \cite{\qaffine, Theorem 2.3, Proposition
2.11, Corollary 1.5}. For $\OlpMmn$, 
$\AnQ$,
$\Oqsymp$ and
$\Oqeucl$, see \cite{\specstrat,
Theorems 5.3, 5.5, 5.8, 5.11, 6.6, 6.8}. The cases
$\OlpGLn$ and
$\OlpSLn$ follow directly. For $\Oq(G)$ and $\Oqp(G,L)$, see (3.1) and
(3.2). \qed\enddemo

\head II. An Axiomatic Development of Stratification Properties\endhead

This part of the paper is a variation on a theme developed in joint work
with E\. S\. Letzter \cite{\specstrat}. The goal is to show that, in the
case of a torus $\H$ acting {\it rationally\/} on a noetherian
$k$-al\-ge\-bra $A$, the $\H$-strat\-i\-fi\-ca\-tion of $\spec A$ enjoys 
analogs of all the properties of the stratification presented in Theorem
3.6, under relatively mild hypotheses on $A$ and $\H$. Both developments
(here and in \cite{\specstrat}) include -- at present -- one hypothesis
which requires substantial work to verify in the examples of interest.
In \cite{\specstrat, Section 6}, that hypothesis is the complete
primeness of the $\H$-prime ideals of $A$. Here, we rely on a type of
normal separation in $\Hspec A$. One disadvantage is the ``loss'' of one
example, namely $\OlpMmn$, for which normal separation has not been
verified. (It is conjectured to hold, and is easily checked for the case
$m=n=2$.) However, the development here provides a better fit with the
work of Hodges-Levasseur, Joseph, and Hodges-Levasseur-Toro on $\Oq(G)$
and $\Oqp(G,L)$ in which the stratified picture of $\spec A$ first
arose. For example, the localizations $(A/J)[\E_J^{-1}]$ in the present
development coincide with the ones obtained in \cite{\HLevone, \HLevtwo,
\Jtwo, \Jbook, \HLT} in the appropriate cases, whereas the localizations
used in
\cite{\specstrat, Section 6} are larger. Here we also obtain a few
additional properties, namely that $(A/J)[\E_J^{-1}]$ is an affine
$k$-al\-ge\-bra in appropriate circumstances, and that normal separation in
$\Hspec A$ by $\H$-ei\-gen\-vec\-tors implies normal separation in $\spec A$.

The development of our theme is based on the equivalence between
rational actions of $(k^\times)^r$ and $\ZZ^r$-gradings (see (5.1) for
details). We begin by working out the basic properties of a
stratification relative to graded-prime ideals for $\ZZ^r$-graded rings
satisfying an appropriate normal separation condition. This may have
some independent interest in the context of group-graded rings. Since it
is at this stage where our variation differs in certain aspects from
that in \cite{\specstrat, Section 6}, we provide full details. This is
the content of Section 4. In Section 5, we translate these results into
the context of $(k^\times)^r$-actions. That process is essentially the
same as in \cite{\specstrat, Section 6}.

\head 4. Graded-nor\-mal separation\endhead

Here we build a stratification in the context of rings graded by a free
abelian group, with a normal separation hypothesis. This hypothesis will
allow us to construct certain localizations which are graded-simple
(i.e., have no nontrivial homogeneous ideals). Hence, we begin by
analyzing the prime spectra of some graded-simple rings.

\proclaim{Lemma 4.1} Let $R$ be a graded-simple ring graded by an
abelian group $G$.

{\rm (a)} The center $Z(R)$ is a homogeneous subring of $R$, strongly
graded by the subgroup $G_Z = \{x\in G \mid Z(R)_x \ne 0\}$ of $G$.

{\rm (b)} Every nonzero homogeneous element of $Z(R)$ is invertible.

{\rm (c)} As $Z(R)$-modules, $Z(R)$ is a
direct summand of $R$.

{\rm (d)} Suppose that $G_Z$ is a free abelian group of finite rank. Choose a
basis $\{g_1,\dots ,g_n\}$ for $G_Z$, and choose a nonzero element $z_j\in
Z(R)_{g_j}$ for each $j$. Then $Z(R)= Z(R)_1 [z_1^{\pm1}, \dots,
z_n^{\pm1}]$, a Laurent polynomial ring over the field $Z(R)_1$. 
\endproclaim

\demo{Proof} It is obvious that $Z(R)$ is a homogeneous subring of $R$,
graded by the set $G_Z$. Part (b)
follows from the graded-simplicity of
$R$, and then it is clear that $G_Z$ is a subgroup of $G$ and that $Z(R)$
is strongly
$G_Z$-graded. This proves part (a). Part (d) is routine (see the proof
of \cite{\specstrat, Lemma 6.3(c)}).

Note that $S= \bigoplus_{x\in G_Z} R_x$ is a homogeneous subring of
$R$ containing $Z(R)$. Further, $S$ is a left $S$-module direct summand
of $R$, with complement
$M=\bigoplus_{x\in G\setminus G_Z} R_x$. On the other hand, as in
\cite{\specstrat, Lemma 6.3(c)}, $S$ is a free $Z(R)$-module with a
basis including $1$, whence $Z(R)$ is a $Z(R)$-module direct summand of
$S$. Part (c) follows. \qed\enddemo

\proclaim{Proposition 4.2} Let $R$ be a graded-simple ring graded by an
abelian group $G$. Then there exist bijections between the
sets of ideals of $R$ and $Z(R)$, given by contraction and
extension, that is,
$$I\mapsto I\cap Z(R) \qquad\qquad \text{and} \qquad\qquad
J\mapsto JR.$$\endproclaim

\demo{Proof} It is clear from Lemma 4.1(c) that $JR\cap Z(R)=J$ for
every ideal $J$ of $Z(R)$. It remains to show that $(I\cap Z(R))R=I$ for
any ideal $I$ of $R$. Set
$J=I\cap Z(R)$, and suppose that $I\ne JR$. Pick an element $x\in
I\setminus JR$ of minimal length, say length $n$. Then $x=x_1+\dots+x_n$
for some nonzero elements $x_i\in R_{g_i}$, where $g_1,\dots,g_n$ are
distinct elements of $G$.

Now $Rx_1R$ is a nonzero homogeneous ideal of
$R$, so $Rx_1R=R$ by graded-simplicity and $\sum_j a_jx_1b_j=1$ for
some $a_j,b_j\in R$. Express each $a_j,b_j$ as a sum of homogeneous
elements, and substitute these expressions in the terms $a_jx_1b_j$. This
expands the sum $\sum_j a_jx_1b_j$ in the form $\sum_t c_tx_1d_t$ with
all $c_t,d_t$ homogeneous. Hence, after relabelling, we may assume that
the
$a_j$ and
$b_j$ are all homogeneous, say of degrees
$e_j$ and $f_j$, respectively.

Comparing identity components in the equation $\sum_j a_jx_1b_j=1$, we
obtain
$$\sum_{e_jg_1f_j=1} a_jx_1b_j=1.$$
Thus, after deleting all other $a_jx_1b_j$ terms, we may assume that
$e_jg_1f_j=1$ for all $j$. Since $G$ is abelian, it follows that $e_jf_j=
g_1^{-1}$ for all $j$. 

Set  $x'= \sum_j a_jxb_j$, and note that $x'\in I$. Further,
$$x'= \sum_{i=1}^n \sum_j a_jx_ib_j$$
where $\sum_j a_jx_ib_j \in R_{g_1^{-1}g_i}$ for each $i$. As a result,
$x'$ is an element whose support is contained in $\{1, g_1^{-1}g_2,
\dots, g_1^{-1}g_n\}$ and whose identity component is $1$. Hence,
$x_1x'$ is an element of $I$ whose support is contained in $\{g_1,\dots,
g_n\}$ and whose $g_1$-component is $x_1$. Comparing this element
with $x$, we see that $x-x_1x'$ is an element of $I$ whose support is
contained in $\{g_2,\dots,g_n\}$. By the minimality of $n$, we must have
$x-x_1x'\in JR$, and so $x'\notin JR$. Therefore, after replacing $x$ by
$x'$, there is no loss of generality in assuming that $g_1=1$ and
$x_1=1$.

\underbar{Claim}: There do not exist $g\in G$ and a nonzero element
$y\in I$ whose support is properly contained in the set
$\{gg_1,\dots,gg_n\}$.

Suppose there do exist such $g$ and $y$. Say the $gg_s$-component of $y$
is nonzero. As above, there exists an element of the form
$y'= \sum_j a_jyb_j$ whose support is properly contained in
$\{g_s^{-1}g_1,
\dots, g_s^{-1}g_n\}$ and whose identity component is $1$. Then
$y'\in I$, and $x-x_sy'$ is an element of $I$ whose support is properly
contained in 
$\{g_1,\dots,g_n\}$.
Since $y'$ and $x-x_sy'$ are elements of $I$ of length less than $n$,
they must lie in $JR$, by the minimality of $n$. But then $x\in JR$,
contradicting our assumptions. Therefore the claim is proved.

Finally, consider any homogeneous element $r\in R$, say $r\in R_g$. Then
$rx-xr$ is an element of $I$ with support contained in $\{gg_2,\dots,
gg_n\}$, and so $rx-xr=0$ by the claim. It follows that $x\in Z(R)$ and
so $x\in J$, contradicting our assumption that $x\notin JR$.

Therefore $I=JR$. \qed\enddemo

\proclaim{Corollary 4.3} If $R$ is a graded-simple ring graded by an
abelian group, then contraction and extension provide mutually inverse
homeomorphisms between $\spec R$ and $\spec Z(R)$. \qed\endproclaim

\definition{4.4. Graded-nor\-mal separation} Let $R$ be a group-graded
ring. Recall that a {\it graded-prime} ideal of $R$ is any proper
homogeneous ideal $P$ such that whenever $I,J$ are homogeneous ideals of
$R$ with
$IJ\subseteq P$, then either $I\subseteq P$ or $J\subseteq P$. We say
that
$R$ has {\it graded-nor\-mal separation} provided that for any proper
inclusion
$P\supset Q$ of graded-prime ideals of $R$, there exists a homogeneous
element $c\in P\setminus Q$ which is normal modulo $Q$.
\enddefinition

\proclaim{Theorem 4.5} Let $R$ be a right noetherian ring graded by
an abelian group
$G$, and assume that $R$ has graded-nor\-mal separation.

Let $J$ be a graded-prime ideal of $R$, set
$$\SS_J= \{P\in\spec R\mid J \text{\ is the largest homogeneous ideal
contained in\ } P\},$$
and let $\E_J$ be the
multiplicative set of all the nonzero homogeneous normal elements in the
ring $R/J$.

{\rm (a)} $\E_J$ is a right and left denominator set of regular elements
in
$R/J$.

{\rm (b)} The localization map
$R\rightarrow R/J\rightarrow R_J= (R/J)[\E_J^{-1}]$
induces a homeomorphism of $\SS_J$ onto $\spec R_J$.

{\rm (c)} Contraction and extension induce mutually inverse
homeomorphisms between $\spec R_J$ and $\spec Z(R_J)$.

{\rm (d)} If $G$ is free abelian of rank $r<\infty$, then $Z(R_J)$ is a
commutative Laurent polynomial ring over the field $Z(R_J)_1$ (the identity
component of $Z(R_J)$ in the induced $G$-grading), in
$r$ or fewer indeterminates.
\endproclaim

\demo{Proof} Without loss of generality, $J=0$. Consequently, $R$ is now
a graded-prime ring, i.e., $0$ is a graded-prime ideal of $R$.

(a) For $x\in \E_J$, observe that
$\rann_R(x)$ is a homogeneous ideal such that $(xR)\rann_R(x) =0$. Since
$R$ is graded-prime and $x\ne 0$, we must have $\rann_R(x)=0$.
Similarly, $\lann_R(x)$=0; thus all elements of $\E_J$ are regular.
The Ore condition follows directly from  normality. 

(b) If $P\in \spec R$ and $P_0$ is the largest homogeneous
ideal contained in $P$, then $P_0$ is a graded-prime ideal. Because of
graded-nor\-mal separation, $P_0$ is nonzero precisely when $P_0\cap \E_J$
is nonempty. Hence,
$\SS_J$ consists precisely of those prime ideals of $R$ disjoint from
$\E_J$. Standard localization theory (e.g., \cite{\GoWa, Theorem 9.22})
now yields part (b).

(c) Since the elements of
$\E_J$ are homogeneous, the
$G$-grading on
$R$ extends (uniquely) to a $G$-grading on $R_J$. Because of the
noetherian hypothesis on $R$, two-sided ideals $I$ in $R$ induce
two-sided ideals $IR_J$ in $R_J$ (e.g., \cite{\GoWa, Theorem 9.20(a)}).
The standard arguments concerning contractions of prime ideals (e.g.,
\cite{\GoWa, Theorem 9.20(c)}) now show that any graded-prime ideal of
$R_J$ must contract to a graded-prime ideal of $R$. Since every nonzero
graded-prime ideal of $R$ meets
$\E_J$, we thus  see that $R_J$ has no nonzero graded-prime ideals. On
the other hand, any maximal proper homogeneous ideal of $R_J$ is
graded-prime. Therefore $R_J$ must be graded-simple. Part (c) now
follows from Corollary 4.3.

(d) This follows from Lemma 4.1(d). \qed\enddemo

\proclaim{Corollary 4.6} Let $R$ be a right noetherian ring graded by
an abelian group
$G$. If $R$ has graded-nor\-mal separation, then $R$ also has normal
separation. \endproclaim

\demo{Proof} Let $P\supset Q$ be a proper inclusion of prime ideals of
$R$. Let $P_0$ and $Q_0$ be the largest homogeneous ideals contained in
$P$ and $Q$, respectively, and note that $P_0$ and $Q_0$ are
graded-prime ideals of $R$. Obviously $P_0\supseteq Q_0$.

If $P_0\ne Q_0$, then by assumption there exists a homogeneous element
$c\in P_0\setminus Q_0$ which is normal modulo $Q_0$. In particular,
$c\in P$ and $c$ is normal modulo $Q$. Since $RcR$ is a homogeneous
ideal of $R$, not contained in $Q_0$, we see that $RcR\not\subseteq Q$.
Thus in this case we are done.

Now assume that $P_0=Q_0$. It is harmless to pass to $R/Q_0$,
and hence there is no loss of generality in assuming that $P_0=Q_0=0$.
Consequently, $0$ is now a graded-prime ideal of $R$.

Let $S$ denote the localization $R_0= R[\E_0^{-1}]$, in the notation of
Theorem 4.5. By part (b) of the theorem, $PS\supset QS$ is a proper
inclusion of prime ideals of $S$. Part (c) then implies that there
exists a central element $z\in PS\setminus QS$. Write $z=cx^{-1}$ for
some $c\in P$ and $x\in \E_0$, and note that $c\notin Q$. Since $xR=Rx$,
we have $zxR=zRx$ in $S$. But $z$ commutes with $R$, and so $zxR=Rzx$,
that is, $cR=Rc$. Therefore $c$ is a normal element of $R$, and normal
separation is proved. \qed\enddemo

\head 5. Stratification under rational torus actions\endhead

We now turn to rational actions by tori as automorphisms of noetherian
algebras. Under suitable normal separation hypotheses, the results of
the previous section apply, yielding a picture of the corresponding
stratifications that incorporates the features of the examples discussed
in Section 3.

Throughout this section, we assume that our base field $k$ is infinite.
In the examples of interest, this is automatic due to the presence of
non-roots of unity.

\definition{5.1. Rational torus actions} Let $A$ be a $k$-al\-ge\-bra, and
let $\H$ be a torus over $k$, say $\H= (k^\times)^r$. An action of $\H$
on $A$ by $k$-al\-ge\-bra automorphisms is said to be {\it rational\/}
provided $A$ is a directed union of finite dimensional $\H$-stable
subspaces $V_i$ such that the induced group homomorphisms $\H\rightarrow
GL(V_i)$ are morphisms of algebraic varieties. Let $\Hhat$ denote the
set of {\it rational characters\/} of $\H$, that is, algebraic group
morphisms $\H\rightarrow k^\times$. Then $\Hhat$ is an abelian group
under pointwise multiplication. Since $k$ is infinite, $\Hhat$ is a
lattice of rank $r$, with a basis given by the component projections
$(k^\times)^r \rightarrow k^\times$.

The rationality of an $\H$-action implies that $A$ is spanned by
$\H$-ei\-gen\-vec\-tors \cite{\Nor, Chapter 5, Corollary to Theorem 36}, and
that the eigenvalues of these eigenvectors are rational. This yields a
$k$-vector space decomposition $A= \bigoplus_{x\in\Hhat} A_x$, where
$A_x$ denotes the $x$-ei\-gen\-space of $A$. (Up to this point, everything
holds for a rational action by $k$-linear transformations.) Since $\H$
acts by automorphisms, $A_xA_y\subseteq A_{xy}$ for all $x,y\in\Hhat$.
Thus we have a grading of $A$ by the free abelian group $\Hhat$, such
that the homogeneous elements are the $\H$-ei\-gen\-vec\-tors in $A$.

Conversely, any grading of $A$ by $\ZZ^r$ arises in this fashion, as
noted in \cite{\ReVo, p\. 784}. To see this, identify $\ZZ^r$ with
$\Hhat$, and let $A= \bigoplus_{x\in\Hhat} A_x$ be an $\Hhat$-grading.
Let each $h\in\H$ act on $A$ as the $k$-linear transformation with
eigenvalue $x(h)$ on $A_x$ for all $x$, that is, $h{.}a= x(h)a$ for 
$a\in A_x$. This yields a rational action of $\H$ on $A$ by $k$-al\-ge\-bra
automorphisms, such that each $A_x$ is the $x$-ei\-gen\-space. The reader
may take this as the definition of a rational $\H$-action if so desired.

In case $k$ is algebraically closed and $A$ is noetherian, a remarkable
theorem of Moeglin-Rentschler and Vonessen \cite{\MoRenone, Th\'eor\`eme
2.12(ii)}, \cite{\VonII, Theorem 2.2} states that $\H$ acts transitively
on the $\H$-strata of rational ideals in $A$, that is, each $\H$-stratum
of rational ideals is a single $\H$-orbit. This holds, in fact, for any
algebraic group $\H$, not just tori. If, further, the Dixmier-Moeglin
equivalence holds, then each $\H$-stratum in $\prim A$ is a single
$\H$-orbit. For the very special case where $\H$ is a torus and all
$\H$-prime ideals of $A$ are completely prime, a direct proof of this
transitivity result is given in \cite{\specstrat, Theorem 6.8}. Here we
shall derive such a result using normal separation in place of complete
primeness -- see Theorem 5.5(c).
\enddefinition

\definition{5.2. Normal separation in $\H$-spec} Assume that we have a
torus $\H$ acting rationally on a $k$-al\-ge\-bra $A$ by $k$-al\-ge\-bra
automorphisms. In this situation, we shall assume that $A$ has been
equipped with the natural $\Hhat$-grading as in (5.1). With respect to
this grading, the homogeneous ideals of $A$ are precisely the
$\H$-stable ideals, and so the graded-prime ideals coincide with the
$\H$-prime ideals. We say that $\Hspec A$ has {\it normal
$\H$-sep\-a\-ra\-tion\/} provided that for any proper inclusion $P\supset Q$
of
$\H$-prime ideals, there exists an $\H$-ei\-gen\-vec\-tor $c\in P\setminus Q$
which is normal modulo $Q$. This condition is just graded-nor\-mal
separation with respect to the $\Hhat$-grading as defined in (4.4).

The term `$\H$-nor\-mal separation' we would reserve for a slightly
stronger condition involving `$\H$-nor\-mal' elements. If $Q$ is an
$\H$-stable ideal of $A$, an element $c\in A$ is said to be {\it
$\H$-nor\-mal modulo $Q$\/} in case there exists $h\in\H$ such that
$ca-h(a)c \in Q$ for all $a\in A$. It is easily checked that each of the
homogeneous components of such an element $c$ is also $\H$-nor\-mal modulo
$Q$. Now $\Hspec A$ has {\it $\H$-nor\-mal separation\/} provided that for
any proper inclusion $P\supset Q$ of
$\H$-prime ideals of $A$, there exists an element $c\in
P\setminus Q$ which is $\H$-nor\-mal modulo $Q$. This strengthening of
normal $\H$-sep\-a\-ra\-tion is not needed in our proofs, but it does hold in
almost all the examples discussed in Sections 2 and 3, and is
conjectured to hold in the remaining one, namely $\OlpMn$.

We note also that since the homogeneous ideals of $A$ coincide with the
$\H$-stable ideals, the largest homogeneous ideal of $A$ contained in a
given ideal $I$ is just $(I:\H)$. Hence, the strata of $\spec A$
with respect to graded-prime ideals, as in Theorem 4.5, coincide with
the
$\H$-strata. Thus Theorem 4.5 and Corollary 4.6
yield the following results.
\enddefinition

\proclaim{Theorem 5.3} Let $A$ be a right noetherian $k$-al\-ge\-bra, and
let $\H$ be a torus of rank $r$, acting rationally on $A$ by $k$-al\-ge\-bra
automorphisms. Assume that $\Hspec A$ has normal $\H$-sep\-a\-ra\-tion. Then
$\spec A$ has normal separation.

Now let $J$ be an $\H$-prime ideal of $A$, let $\E_J$ be the
multiplicative set of all nonzero normal $\H$-ei\-gen\-vec\-tors in $A/J$, and
set $A_J= (A/J)[\E_J^{-1}]$.

{\rm (a)} The localization map $A\rightarrow A/J\rightarrow A_J$ induces
a homeomorphism of $\spec_J A$ onto $\spec A_J$.

{\rm (b)} Contraction and extension induce mutually inverse
homeomorphisms between $\spec A_J$ and $\spec Z(A_J)$.

{\rm (c)} The ring $Z(A_J)$ is a commutative Laurent polynomial ring over the
field
$Z(\fract A/J)^{\H}$ (the fixed subfield of $Z(\fract A/J)$ under the induced
action of $\H$), in
$r$ or fewer indeterminates.
\endproclaim

\demo{Proof} These statements follow immediately from Theorem 4.5 and
Corollary 4.6 except for the description of the coefficient field in part
(c). According to Theorem 4.5(d), this field equals $Z(A_J)_1$. By definition
of the $\Hhat$-grading, $Z(A_J)_1= Z(A_J)^\H$, which is clearly contained in
$Z(\fract A/J)^{\H}$. Thus, it only remains to prove the reverse inequality.

As in the proof of Theorem 4.5(c), $A_J$ is graded-simple with respect to the
$\Hhat$-grading, and so it is $\H$-simple. Given any $u\in Z(\fract
A/J)^{\H}$, observe that the set $I= \{a\in A_J\mid au\in A_J\}$ is a nonzero
$\H$-stable ideal of $A_J$. By $\H$-simplicity, $I=A_J$, whence $u\in A_J$.
Therefore $u\in Z(A_J)^\H$, as desired. \qed\enddemo

\definition{5.4} To compare Theorem 5.3 with \cite{\specstrat, Theorem
6.6}, note that the latter theorem only applies to completely prime
$\H$-prime ideals of $A$. That restriction is due to the lack of any
known graded Goldie theorem for graded-prime rings (cf\.
\cite{\specstrat, \S6.1}), whereas one can easily prove a graded Ore
theorem \cite{\specstrat, Lemma 6.2}. We have finessed the graded Goldie
problem here by assuming a suitable supply of normal elements. Other
advantages of this assumption include the following result.
\enddefinition

\proclaim{Proposition} Let $A$, $\H$, $J$, $\E_J$, $A_J$ be as in
Theorem {\rm 5.3}. If $\Hspec A$ is finite and $A$ is an affine
$k$-al\-ge\-bra, then $A_J$ is an affine $k$-al\-ge\-bra.\endproclaim

\demo{Proof} After passing to $A/J$, we may assume that $J=0$. If $\E_J$
consists of units, then $A_J=A$ and we are done. Otherwise, $A$ contains
some proper nonzero $\H$-stable ideals (generated by nonzero normal
$\H$-ei\-gen\-vec\-tors which are not units). Hence, the maximal proper
$\H$-stable ideals of $A$ are nonzero, and these are, in particular,
$\H$-prime. In other words, 0 is not the only $\H$-prime ideal of $A$.

Let $J_1,\dots,J_n$ be the nonzero $\H$-prime ideals of $A$. By our
normal $\H$-sep\-a\-ra\-tion assumption, each $J_i$ contains a nonzero normal
$\H$-ei\-gen\-vec\-tor $c_i$. Set $c=c_1c_2 \cdots c_n$, a normal
$\H$-ei\-gen\-vec\-tor contained in all nonzero $\H$-primes of $A$. Since
$J=0$ is $\H$-prime, $c\ne0$.

The action of $\H$ on $A$ extends naturally to an action by $k$-al\-ge\-bra
automorphisms on the localization $B=A[c^{-1}]$. Since any $\H$-prime of
$A[c^{-1}]$ contracts to an $\H$-prime of $A$, we see that 0 is the only
$\H$-prime of $B$. Thus $B$ is $\H$-simple.

If $e\in\E_J$, then since $eA$ is an ideal of $A$ and $A$ is right
noetherian, $eB$ is an ideal of $B$. But $eB$ is nonzero and
$\H$-stable, whence $eB=B$. Thus $e$ is right invertible in $B$, and
hence invertible. Therefore $A_J=B$, which is clearly affine.
\qed\enddemo

\definition{5.5} Recall that a noetherian $k$-al\-ge\-bra $A$ satisfies the
{\it Nullstellensatz (over $k$)\/} \cite{\McRob, Chapter 9} provided
that $A$ is a Jacobson ring and the endomorphism rings of all simple
$A$-modules are algebraic over $k$. In that case, it follows from the
Jacobson condition that every locally closed prime of $A$ is primitive,
and from
\cite{\Dix, Lemma 4.1.6} that every primitive ideal of $A$ is rational.
As in \cite{\specstrat, Section 6}, the presence of the Nullstellensatz
yields the following refinements to Theorem 5.3.
\enddefinition

\proclaim{Theorem} Let $A$ be a noetherian $k$-al\-ge\-bra satisfying the
Nullstellensatz, and let $\H$ be a torus  acting rationally on $A$ by
$k$-al\-ge\-bra automorphisms. Assume that $\Hspec A$ is finite and satisfies
normal $\H$-sep\-a\-ra\-tion.

{\rm (a)} $A$ satisfies the Dixmier-Moeglin equivalence, and the
primitive ideals of $A$ are precisely the maximal elements of the
$\H$-strata in $\spec A$.

{\rm (b)} The fields $Z(\fract A/J)^{\H}$ occurring in Theorem {\rm 5.3(c)\/}
are all algebraic over $k$.

{\rm (c)} If $k$ is algebraically closed, each $\H$-stratum in $\prim A$
consists of a single $\H$-orbit. \endproclaim

\demo{Proof} Except for minor modifications, the proofs are the same as
in \cite{\specstrat, Theorem 6.8, Corollary 6.9}. 

Since $\Hspec A$ is finite, maximal elements of $\H$-strata are locally
closed in $\spec A$ \cite{\specstrat, \S2.2(ii)}. Hence, to prove part
(a) it remains to show that rational ideals are maximal within their
$\H$-strata in $\spec A$.

Let $J\in\Hspec A$, and define $\E_J$ and $A_J$ as in Theorem 5.3.
Then $Z(A_J)$ is a Laurent polynomial ring of the form $k_J[z_1^{\pm1},
\dots, z_n^{\pm1}]$, where $k_J= Z(\fract A/J)^{\H}$. By the previous
paragraph, any maximal element of
$\spec_J A$ is locally closed, and such primes are rational because of
the Nullstellensatz. Hence, there exist rational ideals in $\spec_JA$.

Now let $P$ be any rational ideal in $\spec_JA$, and set $Q= PA_J\cap
Z(A_J)$. Then $Z(A_J)/Q$ embeds in $Z(\fract A/P)$, which is algebraic
over $k$ by assumption. Hence, $k_J$ is algebraic over $k$ (thus
establishing part (b)), and $Q$ is a maximal ideal of $Z(A_J)$. It
follows that $P$ is maximal in $\spec_JA$, which completes part (a).

If $k$ is algebraically closed, one checks that the induced action of
$\H$ on $Z(A_J)$ incorporates all automorphisms such that $z_1\mapsto
\beta_1z_1$, \dots, $z_n\mapsto \beta_nz_n$ for arbitrary $\beta_i\in
k^\times$. Hence, $\H$ acts transitively on $\max Z(A_J)$, and part (c)
follows. \qed\enddemo

\definition{5.6} Theorems 5.3 and 5.5 can be used to derive all cases of
Theorem 3.6 except that of $\OlpMmn$. All the necessary hypotheses except
for normal separation can be verified fairly readily from the basic
descriptions of the algebras. Further, these hypotheses are also known
to hold in $\OlpMmn$, where  normal separation is conjectured. Thus, we
conclude by focusing on normal separation as a key problem:
\enddefinition

\proclaim{Problem} Find general hypotheses on a noetherian $k$-al\-ge\-bra
$A$ equipped with a rational action of a torus $\H$ by $k$-al\-ge\-bra
automorphisms which

{\rm (a)} imply that $\Hspec A$ is finite and satisfies normal
$\H$-sep\-a\-ra\-tion;

{\rm (b)} are readily proved for the algebras described in Section 1.
\endproclaim

\head Acknowledgement\endhead

We thank Ken Brown, Iain Gordon, Tim Hodges, David Jordan and Thierry
Levasseur for helpful correspondence and references.

\enddocument